\documentclass[a4paper,reqno]{amsart}

\usepackage{amsmath}
\usepackage{amsthm}
\usepackage{amsfonts}
\usepackage{amssymb}

\usepackage{color}

\usepackage{hyperref}

\hypersetup{
    colorlinks=false,
    pdfborder={0 0 0},
    pdfstartview={XYZ null null 1.00}
}

\numberwithin{equation}{section}

  \def\<{\langle}
  \def\>{\rangle}

  \def\ker{\mathrm{Ker}\,}
  
  \def\inv{\mathrm{Inv}\,}
  \def\inter{\mathrm{int}\,}
  \def\re{\Re\,}
  \def\ve{\varepsilon}

  \newcommand{\x}{\bar u}
  \newcommand{\y}{\bar v}
  \newcommand{\z}{\bar w}

  \def\d{ \, d }

  \def\R{\mathbb{R}}
  
  \def\Q{\mathbb{Q}}
  \def\C{\mathbb{C}}

  \def\ffor{ \qquad \mathrm{for} \quad }
  \def\aas{ \qquad \mathrm{as} \quad }

  \def\A{{\bf A}}
  
  \def\E{{\bf E}}

  \def\P{{\bf P}}
  \def\Q{{\bf Q}}
  
  \def\F{{\bf F}}
  \def\G{{\bf G}}

  \def\Pphi{{\bf \Phi}}
  \def\Ppsi{{\bf \Psi}}

  \def\o{\overline}

  \def\h{\widehat}

\theoremstyle{plain}
  \newtheorem{theorem}{Theorem}[section]

  \newtheorem{proposition}[theorem]{Proposition}
  
  \newtheorem{lemma}[theorem]{Lemma}

\theoremstyle{definition}
  
  \newtheorem{example}[theorem]{Example}
  \newtheorem{remark}[theorem]{Remark}

\begin{document}

\title[Homotopy invariants methods in global...]{Homotopy invariants methods in the global dynamics of strongly damped wave equation}

\author{Piotr Kokocki}
\address{\noindent Faculty of Mathematics and Computer Science \newline Nicolaus Copernicus University \newline Chopina 12/18, 87-100 Toru\'n, Poland}
\address{\noindent BCAM - Basque Center for Applied Mathematics \newline Alameda de Mazarredo 14, 48009, Bilbao, Spain}
\email{pkokocki@mat.umk.pl}
\thanks{The researches supported by the NCN Grant no. 2011/01/N/ST1/05245}



\begin{abstract}
We are interested in the following differential equation
$\ddot u(t) = -A u(t) - c A \dot u(t) + \lambda u(t) + F(u(t))$ where $c > 0$ is a damping factor, $A$ is a sectorial operator and $F$ is a continuous map. We consider the situation where the equation is at resonance at infinity, which means that $\lambda$ is an eigenvalue of $A$ and $F$ is a bounded map. We provide geometrical conditions for the nonlinearity $F$ and determine the Conley index of the set $K_\infty$, that is the union of the bounded orbits of this equation.
\end{abstract}

\maketitle

\setcounter{tocdepth}{2}

\section{Introduction}

We consider the following {\em strongly damped wave equation}
\begin{equation}\label{eq2wave}
\left\{\begin{aligned}
    & u_{tt} = \Delta u + c \Delta u_t + \lambda u + f(x,u), && t\ge 0, \ x\in\Omega \\
    & u(t,x) = 0, && t\ge 0, \ x\in\partial\Omega
\end{aligned}\right.
\end{equation}
where $c> 0$ is a damping factor, $\lambda$ is a real number and $f:\Omega\times\R\to\R$ is a continuous map defined on an open bounded set $\Omega\subset\R^n$. We are interested in the dynamics of the equation \eqref{eq2wave} depending on the parameter $\lambda\in\R$, with the assumption that $f$ is a bounded map. Let us denote by $(\lambda_i)_{i \ge 1}$ the positive sequence of eigenvalues of the operator $A_2$ given by
$$A_2 u := -\Delta u \text{ \ for \ } u\in D(A_2) := H^2(\Omega)\cap H^1_0(\Omega).$$
If $\lambda < \lambda_1$ then the results from \cite{MR1972247}, \cite{MR1778284} say that the equation \eqref{eq2wave}, considered on the space $H^1_0(\Omega)\times L^2(\Omega)$, possess a global attractor i.e. maximal compact invariant set $K$ with the property that after sufficiently large time the trajectory of any bounded set is arbitrary close to $K$. For more details and properties on global attractors see e.g. \cite{MR1778284}, \cite{MR0941371}, \cite{MR1914080}.
If we assume that $\lambda \ge \lambda_1$, then the situation is different. In particular the case when $\lambda\in(\lambda_k,\lambda_{k+1})$ for some $k\ge 1$, was considered in \cite{MR731150}. It was proved that the set $K_\infty$ defined as the union of all bounded orbits of the equation \eqref{eq2wave}, is nonempty, compact and its Conley index is equal to $\Sigma^{d_k}$ -- the homotopy type of the suspension of pointed $d_k$-dimensional sphere, where $d_k:=\sum_{i=1}^k \dim\ker(\lambda_i I -A)$. For more details on Conley index we refer the reader to \cite{MR511133}, \cite{MR637695}, \cite{MR910097}, \cite{MR688146}. In this case $K_\infty$ is a maximal compact invariant set, but it is no longer a global attractor. To be more precise, up to some admissible homotopy, the dynamics of the equation \eqref{eq2wave} in a neighborhood of $K_\infty$ is similar to the dynamics of an isolated invariant set with $d_k$-dimensional unstable manifold.

If we assume that the equation \eqref{eq2wave} is at {\em resonance at infinity}, that is,
$$\ker(\lambda I - A_2)\neq \{0\} \text{ \ and \ } f \text{ \ is a bounded map,}$$
then the problem of existence of a nonempty maximal compact invariant set seems to be not considered so far.
In this case the situation becomes more complicated because there are examples of the nonlinearity $f$ such that the maximal compact invariant set does not exist and, what is worse, the equation \eqref{eq2wave} does not admit any bounded orbit. This fact will be explained in detail in Remark \ref{rem-non-ex3}.
Therefore, the goal of this paper is to study topological properties of the set $K_\infty$ by means of the homotopy invariants methods. In particular we introduce geometrical conditions for the nonlinearity $f$ that, in a consequence, will guarantee that the set $K_\infty$ is compact and nonempty.

To this end we will consider the following more general abstract equation
\begin{align}\label{row-par}
\ddot u(t) = -A u(t) - c A \dot u(t) + \lambda u(t) + F(u(t)), \quad t > 0
\end{align}
where $c> 0$ is still a damping factor, $\lambda$ is a real number, $A: X\supset D(A)\to X$ is a sectorial operator on a Banach space $X$ with compact resolvents and $F:X^\alpha\to X$ is a continuous bounded map, where $X^\alpha :=  D(A^\alpha)$ for $\alpha\in(0,1)$, is a fractional space endowed with the graph norm. For more details on construction and properties of fractional spaces we refer the reader to \cite{MR1721989}, \cite{MR610244}, \cite{H-F}, \cite{Pazy}.
After passing into the abstract framework, we will say that the equation \eqref{row-par} is at {\em the resonance at infinity}, provided
$$\ker(\lambda I - A)\neq \{0\} \text{ \ and \ } F \text{ \ is a bounded map.}$$
The resonance phenomenon for the equation \eqref{row-par} is of importance from the point of view of mathematics, physics and engineering. See for example \cite{MR1084570}, \cite{MR1389086} for an extensive discussion on the meaning of resonance in the oscillations of suspension bridges. To explain the methods that we will use in this paper more precisely, observe that the equation \eqref{row-par} can be written in the following form
\begin{equation}\label{row-hyp}
    \dot w(t) = -\A w(t) + \F(w(t)), \qquad t > 0,
\end{equation}
where $\A:\E\supset D(\A)\to\E$ is a linear operator on the space $\E:=X^\alpha\times X$ given by
\begin{equation*}
\begin{aligned}
D(\A) & :=\{(x,y)\in \E \ | \ x + c y\in D(A)\} \\
\A(x,y) & :=(-y,A(x + c y) - \lambda x) \ffor (x,y)\in D(\A),
\end{aligned}
\end{equation*}
and $\F: \E\to\E$ is the map defined by $\F(x,y):=(0,F(x))$ for $(x,y)\in\E$. Assume that, for every initial data $(x,y)\in \E$, the equation \eqref{row-hyp} admits a unique mild solution $w(\,\cdot\,; (x,y)):[0,+\infty)\to \E$ starting at $(x,y)$. Then, we can define the associated semiflow $\Pphi:[0,+\infty)\times\E \to \E$ by
\begin{equation*}
\Pphi(t, (x,y)):= w(t;(x,y)) \ffor t\in[0,+\infty), \ (x,y)\in \E.
\end{equation*}
Using the fact that $A$ has compact resolvent, it was proved in \cite{MR731150} that any bounded set $N\subset\E$
is admissible with respect to the semiflow $\Pphi$, that is, if sequences $(t_n)$ in $[0,+\infty)$ and $(z_n)$ in $\E$ are such that $t_n\to +\infty$ as $n\to +\infty$ and $$w([0,t_n]\times\{z_n\})\subset N \ffor n\ge 1,$$ then the set $\{w(t_n; z_n) \ | \ n\ge 1\}$ is relatively compact in $\E$. Therefore, the natural way to study invariant sets for the equation \eqref{row-hyp} is to exploit the Conley index theory developed for arbitrary metric spaces in \cite{MR637695}, \cite{MR910097}. More precisely, we introduce {\em resonant conditions} for the nonlinearity $F$ and apply them to construct an isolating neighborhood $N$ with the property that every bounded orbit of \eqref{row-hyp} lies in the interior of $N$. Then we compute the Conley index of $\inv (N)$ with respect to $\Pphi$ in the terms of the resonant conditions imposed earlier on the nonlinearity. The non-triviality of Conley index will guarantee that the set $K_\infty$ is compact and non-empty.

Roughly speaking, the resonant conditions read as follows. Let $\lambda\in\R$ be given eigenvalue of the operator $A$ such that its the geometric and algebraic multiplicities are equal. Since $A$ has compact resolvents, the eigenvalue $\lambda$ is isolated and there is a subspace $V\subset X$ with the property that $X = \ker(\lambda I - A)\oplus V$ and $\sigma(A_V) = \sigma(A)\setminus\{\lambda\}$, where $A_V$ is the part of the operator $A$ in $V$. Assume that the space $X$ is continuously embedded in a Hilbert space $H$, equipped with a norm $\|\cdot\|_H$ and a scalar product $\<\,\cdot\, , \,\cdot\,\>_H$, such that the operator $A$ has self-adjoint
extension $\widehat{A}:H \supset D(\widehat{A}) \to H$. We say that condition $(G1)$ is satisfied provided
\begin{equation*}
\ \left\{\begin{aligned}
& \text{for any balls } B_1\subset V\cap X^\alpha, \ B_2\subset \ker(\lambda I - A) \text{ there are } R,\rho > 0 \text{ such that } \\
& \<F(x + y), x\>_H > - \<F(x + y), z\>_H + \rho \\
& \text{for } (y,z)\in B_1 \times B_2 \ \text{and} \ x\in \ker(\lambda I - A) \text{ with } \|x\|_H\ge R.
\end{aligned}\right.
\end{equation*}
In the similar way, we say that condition $(G2)$ is satisfied provided
\begin{equation*}
\ \left\{\begin{aligned}
& \text{for any balls } B_1\subset V\cap X^\alpha, \ B_2\subset \ker(\lambda I - A) \text{ there are } R,\rho > 0 \text{ such that } \\
& \<F(x + y), x\>_H < - \<F(x + y), z\>_H - \rho\\
& \text{for } (y,z)\in B_1 \times B_2 \ \text{and} \ x\in  \ker(\lambda I - A) \text{ with } \|x\|_H\ge R.
\end{aligned}\right.
\end{equation*}
An important property of these conditions is the fact that if $F$ is a Nemitskii operator associated with $f:\Omega\times\R\to\R$, then $(G1)$ and $(G2)$ are implicated by well-known Landesman-Lazer conditions introduced in \cite{MR0267269} as well as by strong resonance conditions considered in \cite{MR713209}.

It is also worth to note that the Landesman-Lazer and strong resonance conditions have been widely used in the last years to obtain stationary points and periodic solutions, especially for the heat and wave equations (see for example \cite{MR0487001}, \cite{MR0513090}, \cite{MR2830707}, \cite{hess}, \cite{MR1055536}, \cite{MR597281}). In this paper we will show that these conditions can be also successfully used to prove the existence of maximal compact invariant sets and connecting orbits for the equation \eqref{eq2wave}. For the other results in this direction see \cite{arieta}, where the authors analyze a resonant problem for the heat equation with the nonlinearity located at the boundary and study the behavior of the set $K_\infty$ as a parameter crosses the first Steklov eigenvalue. We also refer the reader to \cite{MR798176} where the topological properties of the set $K_\infty$ are studied for parabolic equation with resonance at zero and with non-resonance assumption at infinity. For the results concerning the existence of bounded solutions for the ordinary differential equations in the resonant situation see \cite{MR1010407}, \cite{MR1678160}, \cite{MR1342038}.

The paper is organized as follows. In Section 2, we will deal with the spectral properties of the operator $\A$. In Proposition \ref{th-spec-dec} we provide a relationship between spectral decomposition of the operators $\A$ and $A$ that will be fundamental in this paper.

Section 3 is devoted to the mild solutions for the equation \eqref{row-par} where the nonlinearity depends additionally from a parameter. We provide the standard facts concerning the existence and uniqueness of mild solutions and we focus on continuity and compactness properties for the associated semiflow.

In Section 4 we provide geometrical assumptions for the nonlinearity $F$ and prove {\em the Conley index formula for invariant sets}, Theorem \ref{th-lan-laz}, that is the main result of this paper. Finally, in Section 5 we provide applications of the obtained abstract results to partial differential equations. First of all we formulate Theorems \ref{lem-est2} and \ref{lem-est3} claiming that if $F$ is the Nemitskii operator associated with the map $f$, then the well known Landesman-Lazer (see \cite{MR0267269}) and strong resonance conditions (see \cite{MR713209}) are actually particular case of $(G1)$ and $(G2)$. Then, we provide the criteria on the existence of nonempty maximal compact invariant sets and the criteria on existence of orbits connecting stationary points for the equation \eqref{eq2wave} in terms of Landesman-Lazer and strong resonance conditions. \\[5pt]
\noindent {\em \bf Notation and terminology.} Let $A:X\supset D(A)\to X$ be a linear operator on a real Banach space $X$ equipped with the norm $\|\cdot\|$. We say that the operator $A$ is sectorial provided there are $\phi\in(0,\pi/2)$,  $M\ge 1$ and $a\in\R$, such that the sector $$S_{a,\phi}:=\{\lambda\in\C \ | \ \phi \le |\mathrm{arg} \, (\lambda - a) \le \pi, \ \lambda\neq a\}$$ is contained in the resolvent set of $A$ and furthermore $$\|(\lambda I - A)^{-1}\| \le M/ |\lambda - a| \ffor \lambda\in S_{a,\phi}.$$
It is well-known that if $A$ is sectorial, then $-A$ is an infinitesimal generator of analytic semigroup which, throughout this paper, will be denoted by $\{S_A(t)\}_{t\ge 0}$. The operator $A$ is called positive if $\re \mu > 0$ for any $\mu\in\sigma(A)$. It can be proved that, if $A$ is positive and sectorial, then given $\alpha \ge 0$ the integral
\begin{equation*}
A^{-\alpha} := \frac{1}{\Gamma(\alpha)}\int_0^\infty t^{\alpha - 1}S_A(t) \d t.
\end{equation*}
is convergent in the uniform operator topology of the space $\mathcal{L}(X)$. Consequently we can define {\em the fractional space} associated with $A$ as the domain of the inverse operator $X^\alpha:= D(A^\alpha)$. The space $X^\alpha$ endowed with the graph norm $\|x\|_\alpha := \|A^\alpha x\|$ is a Banach space, continuously embedded in $X$. For more details on sectorial operators and fractional spaces, we refer the reader to \cite{MR610244}, \cite{H-F}, \cite{Pazy}.

\section{Spectral decomposition}

In this section we will assume that $A:X\supset D(A)\to X$ is a sectorial operator on a real Banach space $X$, equipped with the norm $\|\cdot\|$, such that the following conditions are satisfied: \\[3pt]
\noindent\makebox[9mm][l]{$(A1)$}\parbox[t][][t]{118mm}{the operator $A$ is positive and has compact resolvents,}\\[3pt]
\noindent\makebox[9mm][l]{$(A2)$}\parbox[t][][t]{118mm}{there is a Hilbert space $H$ endowed with a scalar product $\<\,\cdot\,, \,\cdot\,\>_H$ and a norm $\|\,\cdot\,\|_H$ such that $X$ is embedded in $H$ by a continuous injective map $i:X \hookrightarrow H$,}\\[3pt]
\noindent\makebox[9mm][l]{$(A3)$}\parbox[t][][t]{118mm}{there is a self-adjoint operator $\h A:H\supset D(\h A) \to H$ such that $\mathrm{Gr}\,(A)\subset \mathrm{Gr}\,(\h A)$, where the inclusion is understood in the sense of the product map $i\times i$.}

\begin{remark}\label{rem-pom}
Under the assumptions $(A1)-(A3)$, the spectrum $\sigma(A)$ consists of a sequence (possibly finite) of real positive eigenvalues $$0 < \lambda_1 < \lambda_2 < \ldots < \lambda_i < \lambda_{i+1} < \ldots$$ such that $\dim (\lambda_i I - A) < +\infty$ for $i\ge 1$. \\[2pt] \noindent Indeed, the operator $A$ has compact resolvents which implies that there is a complex sequence $(\lambda_i)_{i\ge 1}$ such that
\begin{equation}\label{nini}
\sigma(A) = \sigma_p(A) = \{\lambda_i\}_{i\ge 1} \text{ and } \dim_\C \ker\,(\lambda_i I - A_\C) < +\infty \quad\text{for} \quad i\ge 1,
\end{equation}
where $A_\C$ is the complexification of the operator $A$ (see \cite{MR1345385}, \cite{MR1204883} for more details). Furthermore the sequence is finite or $|\lambda_i|\to +\infty$ as $n\to +\infty$. Since $\lambda_i\in\sigma_p(A)$ we have also that $\lambda_i\in\sigma_p(\h A)$, for $i\ge 1$, as a consequence of $(A3)$. But $\h A$ is a symmetric operator, which implies that $\lambda_i$ are real numbers. Since the operator $A$ is positive, the sequence $(\lambda_i)_{i\ge 1}$ is also positive and $\lambda_i\to +\infty$ as $i\to+\infty$. Hence, in view of the right part of \eqref{nini}, we obtain $\dim (\lambda_i I - A) < +\infty$ for $i\ge 1$. \hfill $\square$
\end{remark}

Let us introduce the space $\E:= X^\alpha\times X$ equipped with the norm
\begin{equation*}
\|(x,y)\|_\E:= \|x\|_\alpha + \|y\| \ffor (x,y)\in \E.
\end{equation*}
We proceed to study the spectral properties of the operator $\A:\E\supset D(\A)\to\E$ given by the formula
\begin{equation}
\begin{aligned}\label{op-hyp}
D(\A) & :=\{(x,y)\in \E = X^\alpha\times X \ | \ x + c y\in D(A)\}, \\
\A(x,y) & :=(-y, A(x + c y) - \lambda x) \ffor (x,y)\in D(\A),
\end{aligned}
\end{equation}
where $\lambda$ is a real number and $c > 0$. Let us now proceed to study the spectral decomposition of the operator $\A$. We will need the following proposition.

\begin{proposition}{\em (\cite[Theorem 2.3]{Kok2})}\label{th:10}
Let us denote by $\lambda := \lambda_k$ the $k$-th eigenvalue of $A$ and let $X_0 := \ker (\lambda I - A)$. Then there are closed subspaces $X_+$, $X_-$ of $X$ such that $X = X_+\oplus X_-\oplus X_0$ and the following assertions hold. \\[-10pt]
\begin{enumerate}
\item[(i)] We have the following inclusions $$X_-\subset D(A), \ \ A(X_-)\subset X_- \text{ \ \ and \ \ } A(X_+\cap D(A)) \subset X_+.$$ Furthermore $X_-$ is a finite dimensional and $X_-= \{0\}$ if $k=1$ and $$X_-=\ker(\lambda_1 I - A)\oplus\ldots\oplus\ker(\lambda_{k-1} I - A) \text{ \ if \ } k\ge 2.$$
\item[(ii)] If $A_+:X_+\supset D(A_+) \to X_+$ and $A_-:X_-\supset D(A_-) \to X_-$ are parts of the operator $A$ in the spaces $X_+$ and $X_-$, respectively, then $$\qquad\qquad\sigma(A_+) = \{\lambda_i \ | \ i\ge k+1 \} \text{ \ and \ } \sigma(A_-) = \{\lambda_i \ | \ 1\le i \le k-1 \} \text{ \ for \ }k\ge 1.$$
\item[(iii)] The spaces $X_0$, $X_-$, $X_+$ are mutually orthogonal in $H$, which means that $\<i (u_l),i(u_m)\>_H = 0$ for $u_l\in X_l$, $u_m\in X_m$ where $l,m\in\{0,-,+\}$, $l\neq m$. \\[-8pt]
\end{enumerate}
\end{proposition}

\begin{remark}\label{rem-proj}
Consider the decomposition $X = X_+\oplus X_-\oplus X_0$ obtained in Proposition \ref{th:10}. Let us denote by $P, Q_\pm:X\to X$ projections given by
\begin{equation}\label{wzz1ab}
    P x = x_0 \ \text{ and } \ Q_\pm x = x_\pm \ffor x\in X,
\end{equation}
where $x = x_+ + x_0 + x_-$ for $x_i\in X_i$ and $i\in \{0,-,+\}$. Since the components are closed in $X$, the projections are continuous maps. Let us denote $Q:= Q_- + Q_+$. Since the inclusion $X^\alpha\subset X$ is continuous, one can decompose $X^\alpha$ on a direct sum of closed spaces
$$X^\alpha = X_0\oplus X^\alpha_-\oplus X^\alpha_+, \text{ \ where \ } X^\alpha_-:= X^\alpha\cap X_-, \ X^\alpha_+:=X^\alpha\cap X_+.$$
Therefore the projections $P$ and $Q_\pm$ can be also considered as continuous maps $P, Q_\pm:X^\alpha\to X^\alpha$ given for any $x\in X^\alpha$ by the formula \eqref{wzz1ab}.
\end{remark}
In the following proposition we proceed to the spectral decomposition of the operator $\A$.
\begin{proposition}{\em (\cite[Theorem 2.6]{kok4})}\label{th-spec-dec}
Let us denote by $\lambda := \lambda_k$ the $k$-th eigenvalue of $A$ and let $\E_0 := \ker(\lambda I - A)\times\ker(\lambda I - A)$. Then there are closed subspaces $\E_+,\E_-$ of $\E$ such that $\E:=\E_-\oplus\E_0\oplus\E_+$ and the following assertions hold.
\begin{enumerate}
\item[(i)] We have the following inclusions $$\E_-\subset D(\A), \ \ \A(\E_-)\subset \E_- \text{ \ \ and \ \ } \A(\E_+\cap D(\A)) \subset \E_+.$$
Furthermore $\dim \E_- = 0$ for $k=1$ and
$$\dim \E_- = \sum_{i=1}^{k-1} \dim\ker(\lambda_i I - A) \ffor k\ge 2.$$
\item[(ii)] If $\A_+:\E_+\supset D(\A_+)\to\E_+$ and $\A_-:\E_-\supset D(\A_-)\to\E_-$ are parts of $\A$ in the spaces $\E_-$ and $\E_+$, respectively, then $$\qquad\quad\sigma(\A_+)\subset \{z\in\C \ | \ \re z > 0\} \text{ \ and \ }\sigma(\A_-) \subset \{z\in\C \ | \ \re z < 0\}.$$
\item[(iii)] If $\P,\Q_-,\Q_+:\E\to\E$ are projections on the components of the decomposition $\E:=\E_-\oplus\E_0\oplus\E_+$ and $\Q := \Q_- + \Q_-$, then
\begin{align}\label{proj}
\P(x,y) = (Px, Py) \ \text{ and } \ \Q(x,y) = (Qx, Qy) \ \text{ for } \ (x,y)\in\E.
\end{align}
\end{enumerate}
\end{proposition}

It is well-known that $\A$ is a sectorial operator (see for example \cite{MR2571574}, \cite{MR702424}) and hence $-\A$ generates a $C_0$ semigroup $\{S_\A(t)\}_{t\ge 0}$ of bounded operators on $\E$. In the following proposition we collect properties of the semigroup that we will use in this paper (see e.g. \cite[Corollary 2.9]{kok4}).
\begin{proposition}\label{cor-dec-2aba}
Let us denote by $\lambda := \lambda_k$ the $k$-th eigenvalue of the operator $A$ and let $\E=\E_-\oplus\E_0\oplus\E_+$ be the decomposition from Proposition \ref{th-spec-dec}.
\begin{itemize}
\item[(i)] We have the following inclusions $S_\A(t) \E_i \subset \E_i$ for $t\ge 0$ and $i\in\{0,-,+\}$.
In particular, for any $t\ge 0$ and $z\in \E$
\begin{equation}\label{ds1-hyp}
S_\A(t)\P z = \P S_\A(t)z \quad\text{and}\quad S_\A(t)\Q_\pm z = \Q_\pm S_\A(t) z.
\end{equation}
\item[(ii)] The $C_0$ semigroup $\{S_\A(t)|_{\E_-}\}_{t\ge 0}$ can be uniquely extended to a $C_0$ group on $\E_-$ and there are constants $M, \delta > 0$ such that
\begin{alignat}{2}\label{nieee1}
\|S_\A(t)z\|_\E & \le M e^{- \delta t} \|z\|_\E && \quad\text{for}\quad z\in \E_+, \ t\ge 0, \\ \label{nieee2}
\|S_\A(t)z\|_\E & \le M e^{ \delta t} \|z\|_\E && \quad\text{for}\quad z\in \E_-, \ t\le 0.
\end{alignat}
\item[(iii)] If $\A_0:\E_0 \to \E_0$ is a part of $\A$ in $\E_0$, then
\begin{equation}
\begin{aligned}\label{dss1}
\A_0(x,y) = (-y, c\lambda y) \text{ \ \ and \ \ } S_{\A_0}(t)(x,y) = S_\A(t)(x,y)
\end{aligned}
\end{equation}
for $(x,y)\in\E_0$ and $t\ge 0$.
\end{itemize}
\end{proposition}

\section{Continuity and compactness properties of semiflows}

Assume that $A: X\supset D(A)\to X$ is a positively defined sectorial operator with compact resolvents on a separable Banach space $X$, equipped with a norm $\|\cdot\|$. Consider the following family of differential equations
\begin{equation}\label{rownn}
    \ddot u(t) = -A u(t) - c A \dot u(t) + \lambda u(t) + F(s, u(t)), \qquad t\ge 0
\end{equation}
where $\lambda$ is a real number, $c > 0$ is a damping factor, $s\in[0,1]$ is a parameter and $F:[0,1]\times X^\alpha\to X$ is a continuous map satisfying the following conditions: \\[5pt]
\noindent\makebox[9mm][l]{$(F1)$}\parbox[t]{118mm}{for every $x\in X^\alpha$ there is an open neighborhood $V\subset X^\alpha$ of $x$ and constant $L > 0$ such that for any $s\in[0,1]$, $x_1,x_2\in V$ $$\|F(s,x_1) - F(s,x_2)\|\le L \|x_1 - x_2\|_\alpha;$$}\\
\noindent\makebox[9mm][l]{$(F2)$}\parbox[t]{118mm}{there is a constant $c_0  > 0$ such that $$\|F(s,x)\| \le c_0(1 + \|x\|_\alpha)\qquad \mbox{ for }\quad  s\in[0,1], \ x\in X^\alpha.$$}\\
\noindent\makebox[9mm][l]{$(F3)$}\parbox[t][][t]{118mm}{for any bounded $V\subset X^\alpha$, the set $F([0,1]\times V)$ is relatively compact in $X$.}\\[5pt]
Let us note that the equation \eqref{rownn} can be written in the following form
\begin{equation}\label{row-abstr}
    \dot w(t) = -\A w(t) + \F(s, w(t)), \qquad t > 0,
\end{equation}
where $\A:\E\supset D(\A)\to\E$ is given by \eqref{op-hyp} and $\F:[0,1]\times \E\to\E$ is defined by
\begin{equation*}
    \F(s, (x,y)):=(0,F(s, x)) \ffor s\in [0,1], \ (x,y)\in\E.
\end{equation*}

Let $J\subset\R$ be an interval and $s\in[0,1]$. We say that a continuous map $w:J \to \E$ is \emph{a mild solution} of equation \eqref{row-abstr}, provided
\begin{equation*}
w(t) = S_\A(t - t')w(t') + \int_{t'}^t S_\A(t - \tau)\F(s,w(\tau)) \d \tau
\end{equation*}
for every $t,t' \in J$, $t' < t$. In the following proposition we collect useful facts concerning existence, continuity and compactness of solutions for the equation \eqref{row-abstr}.
\begin{proposition}\label{th-exist1}
Under the above assumptions the following assertions hold. \\[-12pt]
\begin{itemize}
\item[(a)] For every $s\in[0,1]$ and $(x,y)\in \E$, the equation \eqref{row-abstr} admits a unique mild solution $w(\,\cdot\,;s, (x,y)):[0,+\infty)\to \E$ starting at $(x,y)$. \\[-8pt]
\item[(b)] If sequences $(x_n,y_n)$ in $\E$ and $(s_n)$ in $[0,1]$ are such that $(x_n,y_n)\to (x_0,y_0)$ in $\E$ and $s_n\to s_0$ as $n\to +\infty$, then
\begin{equation*}
w(t;s_n, (x_n,y_n)) \to w(t;s_0, (x_0,y_0)) \aas n\to +\infty,
\end{equation*}
for $t\ge 0$, and this convergence in uniform on bounded subsets of $[0,+\infty)$.
\item[(c)] Let $N\subset\E$ be a bounded subset and let sequences $(s_n)$ in $[0,1]$,  $(t_n)$ in $[0,+\infty)$ and $(z_n)$ in $\E$ are such that $t_n\to +\infty$ as $n\to +\infty$ and $$w([0,t_n]\times\{s_n\}\times\{z_n\})\subset N \ffor n\ge 1.$$ Then the set $\{w(t_n; s_n, z_n) \ | \ n\ge 1\}$ is relatively compact in $\E$.
\end{itemize}
\end{proposition}
\noindent\textbf{Proof.} The proof of point $(a)$ is a consequence of \cite[Theorem 3.3.3]{MR610244} and
\cite[Corollary 3.3.5]{MR610244}. Points $(b)$ and $(c)$ are consequences of \cite[Proposition 4.1]{cw-kok1}. Proof of point $(c)$ can be also found in \cite[Theorem 5.1]{MR910097}. \hfill $\square$

\section{Conley index formula for invariant sets}

In this section we prove {\em the Conley index formula for invariant sets}, which is the main result of this paper. We consider the following differential equation
\begin{equation}\label{eq-hyp}
    \ddot u(t) = -A u(t) - c A \dot u(t) + \lambda u(t) + F(u(t)), \qquad t\ge 0
\end{equation}
where $c>0$ is a damping constant, $\lambda$ is a real number, $A:X\supset D(A) \to X$ is an operator defined on a separable Banach space $X$ with norm $\|\cdot\|$ and $F\colon X^\alpha \to X$ is a continuous map. We are interested in the situation when the equation \eqref{eq-hyp} is at {\em resonance at infinity}, that is,
$$\ker(\lambda I - A)\neq \{0\} \text{ \ and \ } F \text{ \ is a bounded map.}$$
Throughout this section we assume that conditions $(A1)-(A3)$, $(F1)$, $(F3)$ are satisfied and furthermore:\\[2pt]
\noindent\makebox[9mm][l]{$(F4)$}\parbox[t][][t]{118mm}{there is $m > 0$ such that $\|F(x)\|\le m$ for $x\in X^\alpha$.} \\[5pt]
Let us note that the equation \eqref{eq-hyp} can be written in the following form
\begin{equation}\label{rownanie-gl3}
    \dot w(t) = -\A w(t) + \F(w(t)), \qquad t > 0
\end{equation}
where $\A:\E\supset D(\A)\to\E$ is given by \eqref{op-hyp} and $\F:\E\to\E$ is defined by
\begin{equation*}
    \F(x,y):=(0,F(x)) \ \ \text{for} \ \ (x,y)\in\E = X^\alpha\times X.
\end{equation*}
From Proposition \ref{th-exist1} $(a)$ it follows that, for any $(x,y)\in\E$, there is a unique mild solution $w(\,\cdot\,;(x,y)):[0,+\infty)\to \E$ for the equation \eqref{rownanie-gl3}, starting at $(x,y)$. Let $\Pphi:[0,+\infty)\times \E \to \E$ be the semiflow given by $$\Pphi(t, (x,y)) := w(t; (x,y)) \ffor t\in[0,+\infty), \ (x,y)\in \E$$
and let $N\subset\E$ be a closed bounded set. Observe that Proposition \ref{th-exist1} $(b)$ and $(c)$ imply that the semiflow is continuous and the set $N$ is admissible with respect to $\Pphi$, that is, for any sequences $(t_n)$ in $[0,+\infty)$ and $(z_n)$ in $\E$ such that $t_n\to +\infty$ and  $$w([0,t_n]\times\{z_n\})\subset N \ffor n\ge 1,$$ the set $\{w(t_n; z_n) \ | \ n\ge 1\}$ is relatively compact in $\E$.
We recall that the \emph{maximal invariant set} $\inv(N, \Pphi)$ is the set of points $x\in N$, for which there is a map $\sigma:\R \to N$ such that $\sigma(0)=x$ and $$\Pphi(t,\sigma(s)) = \sigma(t+s) \ffor t\ge 0, \ s\in \R.$$
If the set $\inv(N,\Pphi)$ is contained in the interior of $N$, then we call $N$ \emph{isolating neighborhood}.
Furthermore a set $K\subset \E$ is called \emph{isolated invariant set}, provided there is a closed set $N\subset X$, admissible with respect to $\Pphi$, such that $K=\inv(N)\subset\inter N$ (cf. \cite{MR637695}, \cite{MR910097}).

\begin{remark}\label{rem-non-ex3}
If the equation \eqref{eq-hyp} is at resonance at infinity, then there is a nonlinearity $F$ such that the equation
does not admit bounded orbits. \\[5pt]
To see this, let us take $F(x) := y_0$ for $x\in X^\alpha$, where $y_0\in\ker(\lambda I - A)\setminus\{0\}$. If $w:\R\to \E$ is a bounded orbit for \eqref{rownanie-gl3}, then
\begin{equation}\label{equaa}
w(t) = S_\A(t)w(0) + \int_0^t S_\A(t - \tau)(0,y_0)\d\tau\ffor t\in[0,+\infty).
\end{equation}
Consider the direct sum decomposition $\E:=\E_-\oplus\E_0\oplus\E_+$ from Proposition \ref{th-spec-dec} together with the corresponding continuous projections $\P$, $\Q_+$, $\Q_-$.
Let us equip the space $\E_0=\ker(\lambda I - A)\times \ker(\lambda I - A)$ with the following norm
\begin{equation*}
    \|(x,y)\|_{\E_0} := \left(\|x\|^2_H + \|y\|^2_H\right)^{1/2} \ffor (x,y)\in\E_0.
\end{equation*}
Since $\E_0$ is finite dimensional, there are constants $c_1,c_2 > 0$ such that
\begin{equation}\label{eded}
c_1\|(x,y)\|_{\E} \le \|(x,y)\|_{\E_0} \le c_2\|(x,y)\|_\E \ffor (x,y)\in \E_0.
\end{equation}
Acting on the equation \eqref{equaa} by $\P$ and using \eqref{ds1-hyp}, \eqref{dss1} and \eqref{proj}, we infer that
$$\P w(t) = S_{\A_0}(t)\P w(0) + \int_0^t S_{\A_0}(t - \tau)(0, y_0)\d\tau\ffor t\in[0,+\infty).$$
Since $\A_0$ is a bounded operator on a finite dimensional space $\E_0$, it follows that the map $(u_0, v_0):=\P w$ is of class $C^1$ and
\begin{equation*}
\dot u_0(t) = v_0(t), \quad \dot v_0(t) = - c \lambda v_0(t) + y_0, \ffor t \ge 0.
\end{equation*}
From the above equation it follows that
\begin{align*}
\frac{d}{dt} (\<u_0(t), c \lambda y_0\>_H + \<v_0(t), y_0\>_H) = \<v_0(t), c \lambda y_0\>_H + \<- c \lambda v_0(t) + y_0, y_0\>_H = \|y_0\|_H^2
\end{align*}
for $t \ge 0$ and hence
\begin{equation}\label{ro11}
    \<u_0(t), c \lambda y_0\>_H + \<v_0(t), y_0\>_H - \<u_0(0), c \lambda y_0\>_H + \<v_0(0), y_0\>_H= t\|y_0\|_H^2.
\end{equation}
Using \eqref{eded} together with Schwartz inequality, we infer that
\begin{align*}
\<u_0(t), c \lambda y_0\>_H + \<v_0(t), y_0\>_H & \le \|(u_0(t), v_0(t))\|_{\E_0} \|(c\lambda y_0, y_0)\|_{\E_0} \\ & \le c_2\|w(t)\|_\E \|(c\lambda y_0, y_0)\|_{\E_0}.
\end{align*}
Since the orbit $w$ is bounded, it follows that the left side of \eqref{ro11} is bounded for $t\ge 0$. This is a contradiction with $\|y_0\|_H > 0$ and the proof is completed. \hfill $\square$
\end{remark}

To overcome difficulties presented in Remark \ref{rem-non-ex3} we impose additional conditions on the nonlinearity $F$. Recalling that $X^\alpha_+$, $X^\alpha_-$ and $X_0$ are subspaces from Remark \ref{rem-proj}, we equip the spaces $X^\alpha_+\oplus X^\alpha_-$ and $X_0$ in the norms $\|\cdot\|_\alpha$ and $\|\cdot\|_H$, respectively. Then
we introduce the following geometric conditions: \\[1pt]
\begin{equation*}\leqno{(G1)}
\ \left\{\begin{aligned}
& \text{for any balls } B_1\subset X^\alpha_+\oplus X^\alpha_- \text{ and } B_2\subset X_0 \text{ there are } R,\rho > 0 \text{ such that } \\
& \<F(x + y), x\>_H > - \<F(x + y), z\>_H + \rho\\
& \text{for } (y,z)\in B_1 \times B_2 \text{ and } x\in X_0 \text{ with } \|x\|_H\ge R.
\end{aligned}\right.
\end{equation*}
\begin{equation*}\leqno{(G2)}
\ \left\{\begin{aligned}
& \text{for any balls } B_1\subset X^\alpha_+\oplus X^\alpha_- \text{ and } B_2\subset X_0 \text{ there are } R,\rho > 0 \text{ such that } \\
& \<F(x + y), x\>_H < - \<F(x + y), z\>_H - \rho\\
& \text{for } (y,z)\in B_1 \times B_2 \text{ and } x\in X_0 \text{ with } \|x\|_H\ge R.
\end{aligned}\right.
\end{equation*}

We proceed to {\em the Conley index formula for invariant sets} that is the main result of this paper.
\begin{theorem}\label{th-lan-laz}
Assume that $\lambda = \lambda_k$ is the $k$-th eigenvalue of the operator $A$ and let $K_\infty$ be the union of all bounded orbits of the semiflow $\Pphi$.
Then, there is a closed isolating neighborhood $N\subset \E$, admissible with respect to $\Pphi$ such that $K_\infty=\inv(N,\Pphi)$ and the following assertions hold:
\begin{enumerate}
\item[(i)] if $(G1)$ holds, then $h(\Pphi, K_\infty) = \Sigma^{d_k}$, \\[-9pt]
\item[(ii)] if $(G2)$ holds, then $h(\Pphi, K_\infty) = \Sigma^{d_{k-1}}$.
\end{enumerate}
Here $d_l:=\sum_{i=1}^l \dim\ker(\lambda_i I - A)$ for $l\ge 1$ with the exceptional case $d_0:= 0$ and $h$ denotes the Conley index.
\end{theorem}
\begin{remark}
Observe that $d_l$ is a finite number because by Remark \ref{rem-pom}, we have $\sigma(A)=(\lambda_i)_{i\ge 1}$ and  $\dim\ker(\lambda_i I - A) < +\infty$ for $i\ge 1$. \hfill $\square$
\end{remark}

\subsection{Preparation to the proof of Theorem \ref{th-lan-laz}}

In the proof we will use the following family of differential equations
\begin{equation}\label{roww1}
\dot w(t) = - \A w(t) + \G(s,w(t)), \qquad t > 0
\end{equation}
where $\G:[0,1]\times \E \to \E$ is a map defined by $$\G(s,(x,y)):= (0,G(s,x)) = (0, PF(s Q x + P x) + s QF(s Q x + P x))$$
for $s\in[0,1]$ and $(x,y)\in\E$.

\begin{remark}\label{rem-com11}
It is not difficult to see that $G$ satisfies assumption $(F1)-(F3)$. Furthermore there is a constant $m_0 > 0$ such that
\begin{equation}\label{roow}
\|\G(s,z)\|_\E \le m_0 \ffor s\in[0,1], \ z\in\E,
\end{equation}
which is a consequence of assumption $(F4)$. \hfill $\square$
\end{remark}
By Remark \ref{rem-com11} and Proposition \ref{th-exist1} $(a)$, for any $s\in[0,1]$, we can define the semiflow $\Ppsi^s:[0,+\infty)\times X^\alpha \to X^\alpha$ given by $$\Ppsi^s(t, (x,y)):=w(t;s, (x,y)) \ffor t\in[0,+\infty), \ (x,y)\in \E,$$ where $w(\,\cdot\,;s,(x,y)):[0,+\infty)\to \E$ is a unique mild solution of \eqref{roww1} starting at $(x,y)$. Furthermore, Proposition \ref{th-exist1} $(b)$ and $(c)$ imply that the family of semiflows $\{\Ppsi^s\}_{s\in [0,1]}$ is continuous and any bounded subset of $\E$ is admissible with respect to this family. In the following lemma we prove {\em a priori} estimates for the full solutions of the equation \eqref{roww1}.
\begin{lemma}\label{lem-bound}
There is a constant $R > 0$ such that if $w:\R\to \E$ is a full bounded mild solution of \eqref{roww1}, for some $s\in[0,1]$, then
\begin{equation}\label{fs7h2}
\|\Q w(t)\|_\E \le R \ffor t\in\R,
\end{equation}
where $\Q = \Q_+ + \Q_-$ and $\Q_+$, $\Q_-$ are projections from Proposition \ref{th-spec-dec}.
\end{lemma}
\noindent\textbf{Proof.} Since $w$ is a solution of the semiflow $\Ppsi^s$ for some $s\in[0,1]$, one has
$$\Ppsi^s(t-t',w(t')) = w(t) \ffor t,t'\in\R, \ t\ge t',$$
which in turn implies that
\begin{equation}\label{rowwww}
w(t) = S_\A(t - t')w(t') + \int_{t'}^{t} S_\A(t - \tau)\G(s,w(\tau)) \d \tau \ffor t\ge t'.
\end{equation}
Acting on this equation by $\Q_+$ and using \eqref{ds1-hyp}, we infer that
\begin{equation}\label{nier111}
\Q_+w(t) = S_\A(t - t')\Q_+w(t') + \int_{t'}^t S_\A(t - \tau)\Q_+\G(s,w(\tau)) \d \tau \text{ \ for \ } t \ge t'.
\end{equation}
Therefore, by \eqref{nieee1}, one has
\begin{align*}
\|S_\A(t - t')\Q_+w(t')\|_\E & \le M e^{- \delta (t - t')} \, \|\Q_+w(t')\|_\E \\
& \le M e^{- \delta (t - t')} \|\Q_+\|_{L(\E)} \|w(t')\|_\E\ffor t \ge t'.
\end{align*}
Hence the boundedness of $w$ implies that
\begin{equation}\label{granic}
\|S_\A(t - t')\Q_+w(t')\|_\E \to 0 \aas t'\to - \infty.
\end{equation}
Furthermore, by \eqref{nier111}, \eqref{nieee1} and \eqref{roow}, we obtain
\begin{equation*}
\begin{aligned}
\|\Q_+w(t)\|_\E & \le \|S_\A(t - t')\Q_+w(t')\|_\E  + \int_{t'}^t \|S_\A(t - \tau)\Q_+\G(s,w(\tau))\|_\E \d \tau \\
& \le \|S_\A(t - t')\Q_+w(t')\|_\E + M \int_{t'}^t e^{- \delta (t - \tau)} \, \|\Q_+\G(s,w(\tau))\|_\E \d \tau \\
& \le \|S_\A(t - t')\Q_+w(t')\|_\E + m_0 M \|\Q_+\|_{L(\E)}\int_{t'}^t  e^{- \delta (t - \tau)} \d \tau\\
& = \|S_\A(t - t')\Q_+w(t')\|_\E +  m_0M\|\Q_+\|_{L(\E)} (1 - e^{- \delta (t - t')})/\delta,
\end{aligned}
\end{equation*}
and consequently
\begin{equation*}
\|\Q_+w(t)\|_\E \le \|S_\A(t - t')\Q_+w(t')\|_\E + m_0 M \|\Q_+\|_{L(\E)} (1 - e^{- \delta (t - t')})/\delta.
\end{equation*}
Letting $t'\to -\infty$ and using \eqref{granic}, yields
\begin{equation}\label{njnj}
\|\Q_+w(t)\|_\E \le m_0 M \|\Q_+\|_{L(\E)} /\delta \ffor t\in\R.
\end{equation}
On the other hand, acting on \eqref{rowwww} by the operator $\Q_-$ and using \eqref{ds1-hyp}, we obtain
\begin{equation}\label{nierrr}
\Q_-w(t) = S_\A(t - t')\Q_-w(t') + \int_{t'}^t S_\A(t - \tau)\Q_- \G(s,w(\tau)) \d \tau,
\end{equation}
which implies that
\begin{equation}
S_\A(t' - t) \Q_-w(t) = \Q_-w(t') + \int_{t'}^t S_\A(t' - \tau)\Q_-\G(s,w(\tau)) \d \tau \text{ \ for \ }t \ge t'
\end{equation}
because, by Proposition \ref{cor-dec-2aba} $(ii)$, the semigroup $\{S_\A(t)\}_{t\ge 0}$ can be extended on the space $\E_-$ to a $C_0$ group of bounded operators. Let us note that, by \eqref{nieee2}, one has
\begin{align*}
\|S_\A(t' - t) \Q_-w(t)\|_\E & \le  M \, e^{\delta(t' - t)} \|\Q_-w(t)\|_\E  \\
& \le M \, e^{\delta(t' - t)} \|\Q_-\|_{L(\E)} \|w(t)\|_\E \ffor t \ge t',
\end{align*}
and therefore the boundedness of $w$ implies that
\begin{equation}\label{fs11jkjk}
\|S_\A(t' - t) \Q_-w(t)\|_\E \to 0 \aas t\to + \infty.
\end{equation}
Furthermore, using \eqref{nierrr}, \eqref{nieee2} and \eqref{roow} we derive that
\begin{align*}
\|\Q_-w(t')\|_\E & \le \|S_\A(t' - t) \Q_-w(t)\|_\E + \int_{t'}^t \|S_\A(t' - \tau)\Q_-\G(s,w(\tau))\|_\E \d \tau \\
& \le \|S_\A(t' - t) \Q_-w(t)\|_\E + \int_{t'}^t M e^{\delta (t' - \tau)}\|\Q_-\G(s,w(\tau))\|_\E \d \tau \\
& \le \|S_\A(t' - t) \Q_-w(t)\|_\E + \int_{t'}^t m_0 M \|\Q_-\|_{L(\E)} \, e^{\delta (t' - \tau)} \d \tau \\
& = \|S_\A(t' - t) \Q_-w(t)\|_\E + m_0 M \|\Q_-\|_{L(\E)}\left(1 - e^{\delta (t' - t)}\right)/\delta
\end{align*}
and consequently, for $t > t'$, one has
$$\|\Q_-w(t')\|_\E \le \|S_\A(t' - t) \Q_-w(t)\|_\E + m_0 M \|\Q_-\|_{L(\E)}\left(1 - e^{\delta (t' - t)}\right)/\delta.$$
Letting $t\to +\infty$ and using \eqref{fs11jkjk}, yields
\begin{equation}
\|\Q_-w(t')\|_\E \le m_0 M \|\Q_-\|_{L(\E)} /\delta \ffor t'\in\R,
\end{equation}
which together with \eqref{njnj} gives \eqref{fs7h2} and the proof is completed. \hfill $\square$ \\

\subsection{Proof of Theorem \ref{th-lan-laz}}
Let $\{e_1, e_2, \ldots, e_n\}$, where $n=\dim \ker (\lambda I - A)$, be an orthonormal basis of the space $\ker (\lambda I - A)$ equipped with the norm $\|\cdot\|_H$ and the scalar product $\<\,\cdot\, , \,\cdot\,\>_H$. We recall the space $\E_0=\ker (\lambda I - A)\times \ker (\lambda I - A)$ is equipped with the scalar product and norm, given by
\begin{align*}
& \<(x_1,y_1), (x_2,y_2)\>_{\E_0}=\<x_1,x_2\>_H + \<y_1,y_2\>_H && \text{for}\quad (x_1,y_1), (x_2,y_2)\in\E_0, \\
& \|(x,y)\|_{\E_0} = \left(\|x\|_H^2 + \|y\|_H^2\right)^{1/2} && \text{for}\quad (x,y)\in \E_0.
\end{align*}
Let the set $\{f_1, f_2, \ldots, f_{2n}\}$ be the following basis of the space $\E_0$
$$f_i := \begin{cases} ((c\lambda)^2 + 1)^{-1/2} (c\lambda e_i, e_i) & \text{ for } \ 1\le i \le n, \\[5pt]
(0,e_{i - n}) & \text{ for } \ n+1\le i \le 2n \end{cases}$$
and let us define the linear map $W: \E_0 \to \R^n\times\R^n$ as $W(x,y) = (w_1, w_2)$ where
\begin{equation}
\begin{aligned}\label{wzor}
w_1 & := a (c\lambda x_1 + y_1, c\lambda x_2 + y_2, \ldots, c\lambda x_n + y_n) \text{ and }\\
w_2 & := (y_1, y_2, \ldots, y_n)
\end{aligned}
\end{equation}
with $a:=((\lambda c)^2 + 1)^{-1/2}$ and $x_i := \<x,e_i\>_H$, $y_i := \<y,e_i\>_H$ for $i=1,2,\ldots, n$. In view of the fact that the basis $\{e_1, e_2, \ldots, e_n\}$ is orthonormal, one has
\begin{equation}\label{eqq12}
|w_1| = a \|c\lambda x + y\|_H \qquad\text{ and }\qquad |w_2| = \|y\|_H,
\end{equation}
where $|\cdot|$ is the Euclidean norm in $\R^n$. \\[5pt]
\textbf{Step 1.} We proceed to define an isolating neighborhood for the family $\{\Ppsi^s\}_{s\in[0,1]}$. To this end, from Lemma \ref{lem-bound}, we obtain a constant $R_1>0$ with the property that, if $w=(u,v):\R\to\E$ is a bounded full solution of $\Ppsi^s : \E \to \E$ for $s\in[0,1]$, then
\begin{equation}\label{eqq1}
    \|\Q w(t)\|_\E \le R_1 \ffor t\in\R.
\end{equation}
In view of the fact that $X_0 = \ker(\lambda I - A)$ is a finite dimensional space the norms $\|\cdot\|$ and $\|\cdot\|_H$ are equivalent on this space. Hence, by assumption $(F4)$, there is a constant $m_1>0$ such that
\begin{equation}\label{nier-f-haa}
    \|PF(x)\|_H \le m_1 \ffor x\in X^\alpha.
\end{equation}
Choose $R_2 > 0$ such that
\begin{equation}\label{eqq11}
-c\lambda R^2 + m_1 R < -c\lambda R^2_2 + m_1 R_2 < 0 \ffor R\ge R_2,
\end{equation}
and define the following sets
\begin{align*}
B_1 & := \{y\in X^\alpha_+\oplus X^\alpha_- \ | \ \|y\|_\alpha \le R_1 + 1\}, \\
B_2 & :=\{z\in \ker(\lambda I - A) \ | \ \|z\|_H \le R_2 /(c\lambda)\}.
\end{align*}
Using geometrical conditions $(G1)$, $(G2)$ and orthogonality from Proposition \ref{th:10} $(iii)$, one can find $R_3,\rho > 0$ such that
\begin{equation}\label{ineq3}
\<PF(x + y), x\>_H > - \<PF(x + y), z\>_H + \rho
\end{equation}
for any $(y, z,x)\in B_1\times B_2\times X_0$ with $\|x\|_H \ge R_3$, if condition $(G1)$ is satisfied and
\begin{equation*}
\<PF(x + y), x\>_H < -\<PF(x + y),z\>_H - \rho
\end{equation*}
for any $(y, z, x)\in B_1\times B_2\times X_0$ with $\|x\|_H \ge R_3$, if the condition $(G2)$ is satisfied. Let us denote
\begin{equation}\label{eeq1}
    R_4 := ac\lambda R_3 + aR_2
\end{equation}
and define the set $N\subset \E$ as $N:=N_1\oplus N_2$, where
\begin{align*}
N_1 & := \{(x,y)\in\E_-\oplus\E_+ \ | \ \|(x,y)\|_\E \le R_1 + 1\}, \\
N_2 & := W^{-1} M \ \text{ where } \ M:= \{(w_1,w_2)\in\R^{2n} \ | \ |w_1|\le R_4, \ |w_2|\le R_2\}.
\end{align*}
\textbf{Step 2.} We show that $N$ is an isolating neighborhood for the family $\{\Ppsi^s\}_{s\in[0,1]}$. To this end we prove that every bounded orbit $w:\R\to\E$ of the semiflow $\Ppsi^s:\E\to\E$, where $s\in[0,1]$, is contained in the interior of the set $N$. The choice of $R_1$ and the boundedness of $w$ in the space $\E$ gives $\Q w(t)\in \inter N_1$ for $t\in\R$. Hence, it is enough to show that $\P w(t) \in \inter N_2$ for $t\in\R$. First observe that from \eqref{proj} we have $\Q w(t) = (Q u(t),Q v(t))$ for $t\in\R$, and hence \eqref{eqq1} implies that $\|Q u(t)\|_\alpha + \|Q v(t)\| \le R_1$. This in turn gives
\begin{equation}\label{ineq2}
\|Q u(t)\|_\alpha \le R_1 \ffor t\in\R.
\end{equation}
Let us define the map $(w_1,w_2):\R\to\R^{2n}$ by the formula $$(w_1(t),w_2(t)):= W\P(w(t)) = W(Pu(t),Pv(t)) \ffor t\in\R.$$ From \eqref{wzor} we obtain
\begin{equation}
\begin{aligned}\label{równaniiia}
w_1(t) & := a (c\lambda u_1(t) + v_1(t), c\lambda u_2(t) + v_2(t), \ldots, c\lambda u_n(t) + v_n(t)), \\
w_2(t) & := (v_1(t), v_2(t), \ldots, v_n(t)) \ffor t\in\R,
\end{aligned}
\end{equation}
where $a:= ((\lambda c)^2 + 1)^{-1/2}$ and $$u_i(t) := \<Pu(t),e_i\>_H, \quad v_i(t) := \<Pv(t),e_i\>_H \ffor i=1,2,\ldots, n.$$ Acting by the operator $\P$ on the equation
$$(u(t), v(t)) = S_\A (t - t')(u(t'),v(t')) + \int_{t'}^t S_\A(t-\tau)\G(s,(u(\tau),v(\tau)))\d \tau \ \ \text{for} \ \  t\ge t'$$
and using \eqref{ds1-hyp}, we infer that
\begin{align*}
\P(u(t), v(t)) = S_\A (t - t')\P(u(t'),v(t')) + \int_{t'}^t S_\A(t-\tau)\P\G(s,(u(\tau),v(\tau)))\d \tau
\end{align*}
for $t \ge t'$ and consequently, by \eqref{proj} and \eqref{dss1}, we have
$$(Pu(t), Pv(t)) = S_{\A_0} (t - t')(Pu(t'),Pv(t')) + \int_{t'}^t S_{\A_0}(t -\tau)(0,P F(s Q u(\tau) + P u(\tau)))\d \tau.$$
Since $\A_0$ is bounded and defined on finite dimensional space $\E_0$, the maps $t\mapsto Pu(t), Pv(t)$ are continuously differentiable on $\R$ and
\begin{equation*}
\left\{ \begin{aligned}
\frac{d}{dt} Pu(t) & = P v(t), & \qquad & t \in\R \\
\frac{d}{dt} Pv(t) & = - c \lambda P v(t) + PF(s Q u(t) + P u(t)),  &\qquad& t \in \R.
\end{aligned}\right.
\end{equation*}
Hence, one has
\begin{equation}
\begin{aligned}\label{rownaniia1}
c\lambda \dot u_i(t) + \dot v_i(t) & = c\lambda \left\langle\frac{d}{dt}Pu(t),e_i\right\rangle_H + \left\langle\frac{d}{dt}Pv(t),e_i\right\rangle_H \\
& = c\lambda v_i(t) - c\lambda v_i(t) + \<PF(s Q u(t) + P u(t)), e_i\>_H \\
& = \<PF(s Q u(t) + P u(t)), e_i\>_H \ffor t\in\R, \ 1\le i \le n
\end{aligned}
\end{equation}
and
\begin{align}\label{row22}
\dot v_i(t) = \left\langle\frac{d}{dt}Pv(t),e_i\right\rangle_H = - c\lambda v_i(t) + \<PF(s Q u(t) + P u(t)), e_i\>_H
\end{align}
for $t\in\R$, $1\le i \le n$. We show that
\begin{equation}\label{aass}
|w_2(t)| < R_2 \ffor t\in\R.
\end{equation}
Suppose contrary that there is $t_0 \in\R$ such that $|w_2(t_0)| \ge R_2$. By \eqref{równaniiia} and \eqref{row22}, we infer that
\begin{equation}
\begin{aligned}\label{ineq}
\frac{d}{dt}\frac{1}{2}|w_2(t)|^2 & = \dot w_2(t) \cdot w_2(t) = \sum_{i=1}^n \dot v_i(t) \cdot v_i(t) \\
& = \sum_{i=1}^n - c\lambda v^2_i(t) + v_i(t) \<PF(s Q u(t) + P u(t)), e_i\>_H \\
& = - c\lambda |w_2(t)|^2 + \<P F(s Q u(t) + P u(t)), Pv(t)\>_H \ \ \text{for} \ \ t\in\R.
\end{aligned}
\end{equation}
Furthermore, from \eqref{eqq12} and \eqref{nier-f-haa} it follows that
\begin{align*}
\<PF(s Q u(t) + P u(t)), Pv(t)\>_H \le m_1 \|Pv(t)\|_H = m_1 |w_2(t)|,
\end{align*}
and hence, by \eqref{eqq11} and \eqref{ineq}, we deduce that
\begin{equation}
\begin{aligned}\label{rraa}
\frac{d}{dt}\frac{1}{2}|w_2(t)|^2 & = - c\lambda |w_2(t)|^2 + \<PF(s Q u(t) + P u(t)), Pv(t)\>_H \\
& \le - c\lambda |w_2(t)|^2 + m_1 |w_2(t)| \ffor t\in\R.
\end{aligned}
\end{equation}
Let $t_1:=\inf\{t\le t_0 \ | \ |w_2(s)| \ge R_2 \text{ for }s\in [t,t_0]\}$. If $t_1 > -\infty$ then by the continuity $|w_2(t_1)| \ge R_2$ and, by \eqref{rraa} and \eqref{eqq11}, we find that $$\frac{d}{dt}\frac{1}{2}|w_2(t)|^2_{| t=t_1} \le - c\lambda |w_2(t_1)|^2 + m_1 |w_2(t_1)| \le -c\lambda R_2^2 + m_1 R_2 < 0.$$
Hence, there is $\delta > 0$ such that $|w_2(t)| > R_2$ for $t\in(-\delta + t_1, t_1)$ which contradicts the fact that $t_1 > -\infty$ and gives $t_1 = -\infty$.
It follows that $|w_2(t)| \ge R_2$ for $t \le t_0$ and, using \eqref{rraa} again, we obtain
$$\frac{d}{dt}\frac{1}{2}|w_2(t)|^2 \le - c\lambda |w_2(t)|^2 + m_1 |w_2(t)| \le -c\lambda R_2^2 + m_1 R_2 := K < 0 \ffor t\ \le t_0.$$ Integrating this equation we have
$$|w_2(t_0)|^2 + 2(-K)(-t + t_0) \le |w_2(t)|^2 \ffor t\le t_0,$$
which is a contradiction because we assumed that $w$ is bounded and hence $|w_2(t)|$ is bounded for $t\in(-\infty, 0]$ as well. Therefore \eqref{aass} follows as desired. It remains to prove that
\begin{equation}\label{aass2}
|w_1(t)| < R_4 \ffor t\in\R.
\end{equation}
To this end suppose contrary that $|w_1(t_0)| \ge R_4$ for some $t_0\in\R$. Assume also that condition $(G1)$ is satisfied.
Let $t_1:=\sup\{t\ge t_0 \ | \ |w_1(s)| \ge R_4 \text{ for } s\in[t_0,t]\}$. If $t_1 < +\infty$ then, by the continuity, $|w_1(t_1)| \ge R_4$ and, by \eqref{równaniiia} and \eqref{rownaniia1}, we deduce that
\begin{equation*}
\begin{aligned}
\frac{d}{dt}\frac{1}{2}|w_1(t)|^2  & = \dot w_1(t) \cdot w_1(t) = a^2 \sum_{i=1}^n (c\lambda \dot u_i(t) +  \dot v_i(t))(c\lambda u_i(t) +  v_i(t)) \\
& = a^2 \sum_{i=1}^n (c\lambda u_i(t) +  v_i(t))\<PF(s Q u(t) + P u(t)), e_i\>_H \\
& = a^2 \<PF(s Q u(t) + P u(t)), c\lambda Pu(t) + Pv(t)\>_H \ffor t\in\R.
\end{aligned}
\end{equation*}
From \eqref{eqq12} and \eqref{aass}, it follows that $\|Pv(t_1)\|_H = |w_2(t_1)| \le R_2$ and, by \eqref{eeq1},
\begin{equation*}
\begin{aligned}
& ac\lambda R_3 + aR_2 = R_4 \le |w_1(t_1)| = a \|c\lambda P u(t_1) + P v(t_1)\|_H \\
& \qquad \le a c\lambda \|P u(t_1)\|_H + a\|P v(t_1)\|_H \le a c\lambda \|P u(t_1)\|_H + a R_2,
\end{aligned}
\end{equation*}
which in turn implies that
\begin{equation}\label{eeq2}
    \|P u(t_1)\|_H \ge R_3\qquad\text{and}\qquad \|Pv(t_1)\|_H \le R_2.
\end{equation}
Therefore, if $(G1)$ is satisfied, then by \eqref{ineq2} and \eqref{eeq2}, the inequality \eqref{ineq3} gives
\begin{equation}
\begin{aligned}\label{nn11}
\frac{d}{dt}\frac{1}{2}|w_1(t)|^2_{| t=t_1} & = a^2 \<PF(s Q u(t_1) + P u(t_1)), c\lambda Pu(t_1) + Pv(t_1)\>_H \\
& = a^2 c\lambda\<PF(s Q u(t_1) + P u(t_1)), Pu(t_1)\>_H \\
& \qquad +  a^2 c\lambda\<PF(s Q u(t_1) + P u(t_1)),Pv(t_1)/c\lambda\>_H \\
& > a^2 c\lambda \rho > 0.
\end{aligned}
\end{equation}
Consequently, there is $\delta > 0$ such that $|w_1(t)| > R_4$ for $t\in(t_1, t_1 + \delta)$, which contradicts the fact that $t_1 < +\infty$ and gives $t_1 = +\infty$.
This implies that $|w_1(t)| \ge R_4$ for $t \ge t_0$ and similarly as in \eqref{nn11} we obtain
\begin{equation*}
\begin{aligned}
\frac{d}{dt}\frac{1}{2}|w_1(t)|^2 & = a^2 c\lambda\<PF(s Q u(t) + P u(t)), Pu(t)\>_H \\
& \qquad +  a^2 c\lambda\<PF(s Q u(t) + P u(t)),Pv(t)/c\lambda\>_H \\
& > a^2 c\lambda\rho \ \text{ for } \ t\ge t_0.
\end{aligned}
\end{equation*}
Integrating this equation gives
$$|w_1(t)|^2 \ge 2 a^2 c\lambda\rho (t - t_0) + |w_1(t_0)|^2 \ffor t\ge t_0.$$
It is a contradiction because we assumed $w$ is bounded and hence $|w_1(t)|$ is bounded for $t\in[0,+\infty)$ as well. Therefore \eqref{aass2} follows as claimed. Combining \eqref{aass} and \eqref{aass2} yields $\P w(t) \in \inter N_2$ for $t\in\R$ and the set $N$ is isolating neighborhood for $\Ppsi^s$ for any $s\in[0,1]$, provided condition $(G1)$ is satisfied. In particular we proved that the set $K_\infty\subset \inter N$. Similarly, if $(G2)$ holds then the same conclusion can be obtain in the exactly the same way. For brevity of the proof we omit details. \\[5pt]
\textbf{Step 3.}
Let us denote by $\varphi_2:[0,+\infty)\times \E_0\to \E_0$ the semiflow for the equation
\begin{equation}\label{eeq4}
(\dot u(t), \dot v(t)) = - \A_0 (u(t), v(t)) + (0,PF(u(t))) , \qquad t > 0
\end{equation}
and let us define
\begin{align*}
M^i & := \{(w_1,w_2)\in\R^{2n} \ | \ |w_1| < R_4, \ |w_2| = R_2\}, \\
M^e & := \{(w_1,w_2)\in\R^{2n} \ | \ |w_1| = R_4, \ |w_2| < R_2\}, \\
M^b & := \{(w_1,w_2)\in\R^{2n} \ | \ |w_1| = R_4, \ |w_2| = R_2\}.
\end{align*}
Let us denote $B:=N_2$ and let $B^e$, $B^i$ and $B^b$ be the sets of {\em strict egress}, {\em ingress} and {\em bounce off points}, respectively (cf. \cite[Definition 3.2]{MR910097}). We prove that, if condition $(G1)$ is satisfied, then the set $B$ is an isolating block for $\varphi_2$ such that $B^e = W^{-1} M^e$, $B^i = W^{-1} M^i$ and $B^b = W^{-1} M^b$, respectively.

Assume that condition $(G1)$ is satisfied and let $(u,v):[-\delta_2, \delta_1)\to \E_0$ be a solution of the semiflow $\varphi_2$ such that $(u(0),v(0))\in W^{-1} M^i$. Let $(w_1,w_2):\R\to\R^{2n}$ be a map given by the formula
\begin{equation}\label{bhbh3}
(w_1(t),w_2(t)) = W(u(t),v(t)) \ffor t\in [-\delta_2, \delta_1).
\end{equation}
From \eqref{wzor} it follows that
\begin{equation}
\begin{aligned}\label{równaniiia33}
w_1(t) & := a (c\lambda u_1(t) + v_1(t), c\lambda u_2(t) + v_2(t), \ldots, c\lambda u_n(t) + v_n(t)), \\
w_2(t) & := (v_1(t), v_2(t), \ldots, v_n(t)) \ffor t\in [-\delta_2, \delta_1),
\end{aligned}
\end{equation}
where $a:= ((\lambda c)^2 + 1)^{-1/2}$ and $u_i(t) := \<u(t),e_i\>_H$, $v_i(t) := \<v(t),e_i\>_H$ for $i=1,2,\ldots, n$. Since $(u,v)$ is of class $C^1$ it satisfies the equations
\begin{equation}
\left\{ \begin{aligned}
\dot u(t) & = v(t), & \qquad & t \in [-\delta_2, \delta_1), \\
\dot v(t) & = - c \lambda v(t) + PF(u(t)),  &\qquad& t \in [-\delta_2, \delta_1),
\end{aligned}\right.
\end{equation}
which in turn, implies that
\begin{equation}
\begin{aligned}\label{wzorkiii}
c\lambda \dot u_i(t) + \dot v_i(t) & = c\lambda \left\langle\dot u(t),e_i\right\rangle_H + \left\langle \dot v(t),e_i\right\rangle_H \\
& = c\lambda v_i(t) - c\lambda v_i(t) + \<PF(u(t)), e_i\>_H \\
& = \<PF(u(t)), e_i\>_H
\end{aligned}
\end{equation}
for $t\in[-\delta_2, \delta_1)$, $1\le i \le n$ and
\begin{align*}
\dot v_i(t) & =  \left\langle\dot v(t),e_i\right\rangle_H  = - c\lambda v_i(t) + \<PF(u(t)), e_i\>_H
\end{align*}
for $t\in[-\delta_2, \delta_1)$ and $1\le i \le n$. Then $(w_1(0), w_2(0)) = W(u(0),v(0))\in M^i$ and
\begin{equation}
\begin{aligned}\label{sdf}
\frac{d}{dt}\frac{1}{2}|w_2(t)|^2 & = \dot w_2(t) \cdot w_2(t) = \sum_{i=1}^n \dot v_i(t) \cdot v_i(t) \\
& = \sum_{i=1}^n - c\lambda v_i(t)^2 + v_i(t) \<PF(u(t)), e_i\>_H \\
& = - c\lambda \|v(t)\|_H^2 + \<PF(u(t)), v(t)\>_H \ffor t\in [-\delta_2, \delta_1).
\end{aligned}
\end{equation}
Further, by \eqref{eqq12} and \eqref{nier-f-haa}, we have
\begin{align*}
\<PF(u(t)), v(t)\>_H \le m_1 \|v(t)\|_H = m_1 |w_2(t)| \ffor t\in [-\delta_2, \delta_1),
\end{align*}
and hence, by \eqref{sdf} and \eqref{eqq11}, we deduce that
\begin{equation}\label{eeq14}
\begin{aligned}
\frac{d}{dt}\frac{1}{2}|w_2(t)|^2_{| t=0} & = - c\lambda |w_2(0)|^2 + \<PF(u(0)), v(0)\>_H \\
& \le - c\lambda |w_2(0)|^2 + m_1 |w_2(0)| = -c\lambda R_2^2 + m_1 R_2 < 0.
\end{aligned}
\end{equation}
This implies that there are $\ve_1\in (0,\delta_1)$ and $\ve_2 \in (0, \delta_2)$ (when $\delta_2 > 0$) such that $|w_2(t)| < R_2$ for $t\in(0,\ve_1]$ and $|w_2(t)| > R_2$ for $t\in[-\ve_2,0)$. Taking $\ve_1 > 0$ smaller if necessary we also have $w_1(t)\in B_{\R^n}(0,R_4)$, where we define $$B_{\R^n}(0,r):=\{x\in\R^n \ | \ |x| < r\}.$$ Hence $(u(t),v(t))\in\inter B$ for $t\in(0,\ve_1]$ and $(u(t),v(t))\not\in B$ for $t\in[-\ve_2,0)$, which proves that $W^{-1}M^i$ is contained in $B^i$.

Suppose that $(u,v):[-\delta_2, \delta_1)\to \E_0$ is a solution for the semiflow $\varphi_2$ such that $(u(0),v(0)) \in W^{-1} M^e$. Similarly as before, we define a map $(w_1,w_2):[-\delta_2, \delta_1)\to \R^{2n}$ by the formula \eqref{bhbh3}. Then $(w_1(0), w_2(0)) = W(u(0),v(0))\in M^e$, which together with \eqref{równaniiia33} and \eqref{wzorkiii} gives
\begin{align*}
\frac{d}{dt}\frac{1}{2}|w_1(t)|^2  & = \dot w_1(t) \cdot w_1(t) = a^2 \sum_{i=1}^n (c\lambda \dot u_i(t) +  \dot v_i(t))(c\lambda u_i(t) +  v_i(t)) \\
& = a^2 \sum_{i=1}^n (c\lambda u_i(t) +  v_i(t))(PF(u(t)), e_i) \\
& = a^2 (PF(u(t)), c\lambda u(t) + v(t)) \ffor t\in [-\delta_2, \delta_1).
\end{align*}
From \eqref{eqq12} it follows that $\|v(0)\|_H = |w_2(0)| \le R_2$ and
\begin{equation*}
\begin{aligned}
& ac\lambda R_3 + aR_2 = R_4 = |w_1(0)| = a \|c\lambda u(0) + v(0)\|_H \\
& \qquad \le a c\lambda \|u(0)\|_H + a\|v(0)\|_H \le a c\lambda \|u(0)\|_H + a R_2,
\end{aligned}
\end{equation*}
which in turn implies that
\begin{equation}\label{eeq2bb}
    \|u(0)\|_H \ge R_3\qquad\text{and}\qquad \|v(0)\|_H \le R_2.
\end{equation}
Since condition $(G1)$ is satisfied, by \eqref{eeq2bb}, the inequality \eqref{ineq3} yields
\begin{equation}\label{eeq12}
\begin{aligned}
\frac{d}{dt}\frac{1}{2}|w_1(t)|^2_{|t=0} & = a^2 (PF(u(0)), c\lambda u(0) + v(0)) \\
& = a^2 c\lambda(PF(u(0)), u(0)) + a^2c\lambda (PF(u(0)), v(0)/c\lambda) > 0.
\end{aligned}
\end{equation}
From this there are $\ve_1\in (0,\delta_1)$ and $\ve_2 \in (0, \delta_2)$ (when $\delta_2 > 0$) such that $|w_1(t)| > R_4$ for $t\in(0,\ve_1)$ and $|w_1(t)| < R_4$ for $t\in(0,-\ve_2]$. Taking again $\ve_2 > 0$ smaller if necessary, we have also $w_2(t)\in B_{\R^n}(0,R_2)$ for $t\in[-\ve_2, 0)$. Therefore $(u(t),v(t))\in\inter B$ for $t\in[-\ve_2, 0)$ and $(u(t),v(t))\not\in B$ for $t\in(0,\ve_1)$ which implies that the set $W^{-1}M^e$ is contained in $B^e$.

Let $(u,v):[-\delta_2, \delta_1)\to \E_0$ be a solution of the semiflow $\varphi_2$ such that $(u(0),v(0)) \in W^{-1} M^b$.
Let us again define the map $(w_1,w_2):[-\delta_2, \delta_1)\to \R^{2n}$ by \eqref{bhbh3}. Then $(w_1(0), w_2(0)) = W(u(0),v(0))\in M^b$ and the both inequalities \eqref{eeq14} and \eqref{eeq12} hold. Hence there are $\ve_1\in (0,\delta_1)$ and $\ve_2 \in (0, \delta_2)$ such that $|w_1(t)| > R_4$ for $t\in(0,\ve_1)$ and $|w_2(t)| > R_2$ for $t\in[-\ve_2,0)$ (when $\delta_2 > 0$). Therefore  $(u(t),v(t))\not\in B$ for $t\in[-\ve_2, 0)\cup (0,\ve_1)$ and therefore the set $W^{-1}W^b$ is contained in $B^b$.

Since the sets $B^i$, $B^e$ and $B^b$ are mutually disjoint and $\partial B = W^{-1}M^e \cup W^{-1}M^i \cup W^{-1}M^b$ we find that $W^{-1}M^i = B^i$, $W^{-1}M^e = B^e$ and $W^{-1}M^b = B^b$. Consequently $B=N_2$ is an isolating block for the semiflow $\varphi_2$ with $B^- :=B^e\cup B^b= W^{-1} (M^e\cup M^b)$.

Similarly, we can verify that condition $(G2)$ implies that the set $B$ is an isolating block for the semiflow $\varphi_2$ with the boundary $\partial B$ consisting of the strict ingress points. The proof of this fact will be identical with the only difference that in the case of condition $(G2)$ the inequality \eqref{eeq12} will be opposite. \\[5pt]
\textbf{Step 4.} For any $s\in [0,1]$ write $K_s := \inv (\Ppsi^s , N)$. By Step 2 and the homotopy invariance of the Conley index (cf. \cite[Theorem 12.2]{MR910097}) we have
\begin{equation}\label{eeq3}
h(\Pphi, K_\infty) = h(\Ppsi^1, K_1) = h(\Ppsi^0, K_0).
\end{equation}
Let us observe that, the family $\{\Ppsi^s\}_{s\in [0,1]}$ is a homotopy between $\Ppsi^1 = \Pphi$ and $\Ppsi^0$. Furthermore, every solution $(u,v):[0,+\infty) \to \E$ of the semiflow $\Ppsi^0$ satisfies the following formula
$$(u(t), v(t)) = S_\A(t)(u(0), v(0)) + \int_0^t S_\A(t - \tau)(0, PF(P u(\tau)))\d\tau \ffor t\ge 0.$$
Let us denote by $\varphi_1:[0,+\infty)\times \E_+ \oplus \E_- \to \E_+\oplus \E_-$ the semiflow given by $$\varphi_1(t,(x,y)):= S_\A(t) (x,y)\ffor t\in [0,+\infty), \ (x,y)\in \E_+ \oplus \E_-.$$
Then we can easily check that, for any $t\in[0,+\infty)$ and $(x,y)\in \E$, one has $$\Ppsi^0(t, (x,y)) = \varphi_1(t, \Q(x,y)) + \varphi_2(t, \P(x,y)),$$ which implies that $\Ppsi^0$ is topologically equivalent with the cartesian product of $\varphi_1$ and $\varphi_2$. This means that, for any $t\ge 0$ and $(z_1,z_2)\in (\E_-\oplus \E_+) \times \E_0$
\begin{equation}\label{rown1hh}
    \Ppsi^0(t, U(z_1,z_2)) = U(\varphi_1(t, z_1), \varphi_2(t, z_2)),
\end{equation}
where the linear homeomorphism $U:(\E_-\oplus \E_+) \times \E_0 \to \E$ is given by
$$U(z_1,z_2) = z_1 + z_2 \text{ \ for \ } (z_1,z_2)\in (\E_-\oplus \E_+) \times \E_0.$$
Let us denote $K'_1:=\inv (\varphi_1, N_1)$ and $K'_2:= \inv(\varphi_2, N_2)$. In view of \eqref{nieee1}, \eqref{nieee2} and \cite[Theorem 11.1]{MR910097} it follows that $K'_1 = \{0\}$ is an isolated invariant set and
\begin{equation}\label{eeq11}
    h(\varphi_1, K'_1) = \Sigma^{\dim \E_-} = \Sigma^{d_{k-1}},
\end{equation}
where the last equality is a consequence of Proposition \ref{th-spec-dec} $(i)$. Furthermore, by Step 3, we infer that $K'_2$ is an isolated invariant set. Hence the multiplication property of the homotopy index implies that   $K'_1\times K'_2\subset(\E_-\oplus \E_+)\times \E_0$ is an isolated invariant set and
\begin{equation}\label{eeq15}
    h(\varphi_1\times\varphi_2, K'_1\times K'_2) = h(\varphi_1, K'_1)\wedge h(\varphi_2, K'_2).
\end{equation}
Let us denote $K':=\inv (\varphi_1\times\varphi_2, N_1\times N_2)$. Then $K' = K'_1\times K'_2$ and \eqref{rown1hh} implies that $U(K') = \inv (\Ppsi^0, N) = K_0$. Therefore, by the topological invariance of Conley index we find that
\begin{equation}\label{eeq5}
    h(\Ppsi^0, K_0) = h(\varphi_1\times\varphi_2, K'_1\times K'_2).
\end{equation} \\[-5pt]
\textbf{Step 5.} Combining \eqref{eeq3}, \eqref{eeq5}, \eqref{eeq15} and \eqref{eeq11} yields
\begin{equation}\label{eeq5bbb}
    h(\Ppsi, K_\infty) = h(\varphi_1, K'_1)\wedge h(\varphi_2, K'_2) = \Sigma^{d_{k-1}} \wedge h(\varphi_2, K'_2).
\end{equation}
If $(G1)$ holds, then the pair $(B,B^-)$ is homeomorphic with $(M,M^-)$, where
\begin{align*}
M & := \{(w_1,w_2)\in\R^{2n} \ | \ |w_1|\le R_4, \ |w_2|\le R_2\}, \\
M^- & := \{(w_1,w_2)\in\R^{2n} \ | \ |w_1| = R_4, \ |w_2|\le R_2\}
\end{align*}
and therefore
\begin{equation*}
    h(\varphi_2, K'_2) = \Sigma^n = \Sigma^{\dim\ker(\lambda_k I - A)},
\end{equation*}
which together with \eqref{eeq5bbb} gives
\begin{equation}
    h(\Ppsi, K_\infty) = \Sigma^{d_{k-1}}\wedge \Sigma^{\dim\ker(\lambda_k I - A)} = \Sigma^{d_k},
\end{equation}
and the proof of $(i)$ is completed. If condition $(G2)$ is satisfied, then $B:=N_2$ is an isolating block for the semiflow $\varphi_2$ with that boundary $\partial B$ consisting of strict ingress points. In this case the pair $(B,B^-)$ is homeomorphic with $(M, \emptyset)$, which implies that
\begin{equation*}
    h(\varphi_2, K'_2) = \Sigma^0.
\end{equation*}
Combining this with \eqref{eeq5bbb} we find that
\begin{equation}
    h(\Ppsi, K_\infty) = \Sigma^{d_{k-1}}\wedge \Sigma^0 = \Sigma^{d_{k-1}}
\end{equation}
and the point $(ii)$ follows as desired. \hfill $\square$ \\

\section{Applications}

In this section we provide applications of the obtained abstract results to partial differential equations. We will assume that $\Omega\subset\R^n$, $n\ge 1$, is an open bounded set with the boundary $\partial\Omega$ of class $C^1$. Consider the strongly damped wave equation
\begin{equation}\label{A-eps-res-ah}
\left\{\begin{aligned}
& u_{tt}= - c\mathcal{A} \, u_t - \mathcal{A} \, u + \lambda  u + f(x, u), && t > 0, \ x\in\Omega, \\
& u(t,x) = 0 && t \ge 0, \ x\in\partial\Omega,
\end{aligned}\right.
\end{equation}
where $c> 0$ is a damping factor, $\lambda$ is a real number and $\mathcal{A}$ is a differential operator of the following form
$$\mathcal{A} \x (x) = -\sum_{i,j=1}^n D_j(a_{ij}(x)D_i \x(x)) \ffor \x\in C^2(\overline{\Omega}),$$ which is {\em symmetric} $a_{ij} = a_{ji}\in C^1(\overline\Omega)$ and {\em uniformly elliptic} i.e.
$$\sum_{1\le i,j\le n}a_{ij}(x)\xi^i\xi^j \ge c_0 |\xi|^2 \ffor x\in\Omega, \ \xi\in\R^n, \ \ \text{where} \ \ c_0 > 0.$$ Furthermore we assume that $f:\Omega\times\R\to\R$ is a continuous map such that: \\[2pt]
\noindent\makebox[22pt][l]{$(E1)$} \parbox[t][][t]{118mm}{there is $L> 0$ such that if $x\in\Omega$ and $s_1,s_2\in\R$, then
    \begin{align*}
    |f(x,s_1) - f(x,s_2)| \le L |s_1 - s_2|;
    \end{align*}}\\[5pt]
\noindent\makebox[22pt][l]{$(E2)$} \parbox[t][][t]{118mm}{there is $m > 0$ such that $|f(x,s)| \le m$ for $x\in\Omega$, $s\in\R$.}\\[5pt]
Let us introduce the abstract framework for the equation \eqref{A-eps-res-ah}. For this purpose denote $X:=L^p(\Omega)$, for $p\ge 2$, and define the operator $A_p: X\supset D(A_p)\to X$ by
\begin{equation*}
\begin{aligned}
D(A_p) := W^{2,p}(\Omega) \cap W^{1,p}_0(\Omega), \quad A_p \x := \mathcal{A} \x \quad \text{for} \ \ \x\in D(A_p).
\end{aligned}
\end{equation*}
It is known (see e.g. \cite{MR1778284}, \cite{Pazy}, \cite{MR500580}) that $A_p$ is a positive sectorial operator. Let us denote by $X^\alpha := D(A_p^\alpha)$, for $\alpha\in(0,1)$, the associated fractional space and define the map $F\colon X^\alpha \to X$, given for any $\x\in X^\alpha$, by
\begin{equation}
    F(\x)(x) := f(x, \x(x)) \ffor x\in\Omega.
\end{equation}
We call $F$ \emph{the Nemitskii operator} associated with $f$. We are ready to write the equation \eqref{A-eps-res-ah} in the following abstract form
\begin{equation}\label{row-abs}
\ddot u(t)  = - A_p u(t) - A_p \dot u(t) + \lambda u(t) + F (u(t)),  \qquad   t > 0, \\
\end{equation}

\begin{remark}\label{rem-pom2}
$(a)$ We claim that assumptions $(A1)-(A3)$ are satisfied. \\[2pt]
Indeed, $(A1)$ holds because $A_p$ has compact resolvent as it was proved for example in \cite{MR1778284}, \cite{Pazy}, \cite{MR500580}.
To see that $(A2)$ holds it is enough to take $H:=L^2(\Omega)$ equipped with the standard inner product and norm.
Since $p\ge 2$ and $\Omega$ is bounded set, the embedding $i:L^p(\Omega) \hookrightarrow L^2(\Omega)$ is well-defined and continuous. Furthermore, if we define $\h A:= A_2$, then $$i(D(A_p))\subset D(\h A) \quad\text{and}\quad i(A_p \x) = \h A i(\x) \quad\text{for}\quad \x \in D(A_p),$$
which shows that $A_p \subset \h A$ in the sense of the inclusion $i\times i$. Since the operator $\h A$ is self-adjoint (see e.g. \cite{MR1778284}) we see that the assumption $(A3)$ is also satisfied. \\[5pt]
$(b)$ Let us observe that $F$ is satisfies assumptions $(F1)$, $(F3)$ and $(F4)$. \\[2pt]
Indeed, since $f$ satisfies assumptions $(E1)$ and $(E2)$, the fact that conditions $(F1)$, $(F4)$ hold is straightforward. We only show $(F3)$. To this end, take a bounded sequence $(\x_n)$ in $X^\alpha$. Since $A_p$ has compact resolvents, by \cite[Theorem 1.4.8]{MR610244}, the inclusion $X^\alpha\hookrightarrow X$ is compact, and hence, passing if necessary to a subsequence, we can assume that $\x_n \to \x_0$ in $X$ as $n\to +\infty$. Therefore, using $(E2)$ and the dominated convergence theorem, it is not difficult to verify that $F(\x_n) \to F(\x_0)$ in $X$ as $n\to +\infty$, which proves that $(F3)$ holds. \hfill $\square$
\end{remark}

\subsection{Properties of the Nemitskii operator}

We proceed to examine when the the Nemitskii operator $F$ satisfies geometrical conditions $(G1)$ and $(G2)$. Let us first note that, by Remark \ref{rem-pom}, the spectrum $\sigma(A_p)$ consists of sequence of positive eigenvalues $$0 < \lambda_1 < \lambda_2 < \ldots < \lambda_i < \lambda_{i+1} < \ldots$$ which is finite or $\lambda_i \to +\infty$ as $i\to +\infty$. We recall also that $X^\alpha_+$, $X^\alpha_-$ and $X_0$ are subspaces obtained in Remark \ref{rem-proj}, but this time for the operator $A_p$. In particular $X_0 = \ker(\lambda I - A_p)$. Let us start with the following theorem which says that the conditions $(G1)$ and $(G2)$ are implicated by the well-known \emph{Landesman-Lazer} conditions introduced in \cite{MR0267269}.

\begin{theorem}\label{lem-est2}
Let $f_+,f_-\colon \Omega \to \mathbb{R}$ be continuous functions such that
\begin{align*}
f_+(x) = \lim_{s \to +\infty} f(x,s) \quad\text{and}\quad f_-(x) = \lim_{s \to -\infty} f(x,s) \quad\text{for}\quad x\in\Omega.
\end{align*}
\makebox[5mm][l]{(i)} \parbox[t]{117mm}{Assume that
\begin{equation*}
\int_{\{u>0\}} f_+(x) \x(x) \,d x  + \int_{\{u<0\}} f_-(x) \x(x) \,d x > 0 \ \ \text{for} \ \ \x\in X_0\setminus\{0\}.\leqno{(LL1)}
\end{equation*}
If the sets $B_1\subset X^\alpha_+\oplus X^\alpha_-$ and $B_2 \subset X_0$ are bounded in the norms $\|\cdot\|_\alpha$ and $\|\cdot\|_{L^2}$, respectively, then there are constants $R,\rho > 0$ such that
\begin{equation*}
\<F(\z + \x), \x\>_{L^2} > - \<F(\z + \x), \y\>_{L^2} + \rho
\end{equation*}
for any $(\z, \y,\x)\in B_1\times B_2 \times X_0$, with $\|\x\|_{L^2} \ge R$.}\\[5pt]
\makebox[5mm][l]{(ii)} \parbox[t]{117mm}{Assume that
\begin{equation*}
\int_{\{u>0\}} f_+(x) \x(x) \,d x  + \int_{\{u<0\}} f_-(x) \x(x) \,d x < 0 \ \ \text{for} \ \ \x\in X_0\setminus\{0\}.\leqno{(LL2)}
\end{equation*}
If the sets $B_1\subset X^\alpha_+\oplus X^\alpha_-$ and $B_2 \subset X_0$ are bounded in the norms $\|\cdot\|_\alpha$ and $\|\cdot\|_{L^2}$, respectively, then there are constants $R,\rho > 0$ such that
\begin{equation*}
\<F(\z + \x), \x\>_{L^2} < - \<F(\z + \x), \y\>_{L^2} - \rho
\end{equation*}
for any $(\z, \y,\x)\in B_1\times B_2 \times X_0$, with $\|\x\|_{L^2} \ge R$.}\\[5pt]
\end{theorem}
\noindent\textbf{Proof.} Except for technical modifications, the argument goes in the lines of the proof of \cite[Theorem 6.7]{Kok2}. We encourage the reader to reconstruct details. \hfill $\square$ \\

\begin{example}
We describe an example of the situation where the Landesman-Lazer conditions $(LL1)$ and $(LL2)$ are satisfied. \\[5pt]
To this end suppose by the moment that
$$\mathcal{A}\x = -\x_{xx} \quad\text{for}\quad \x\in C^2([0,1]) \quad\text{and}\quad f(s) := \arctan(s) \quad\text{for}\quad s\in\R.$$
Then $A_p u = - u_{xx}$ is defined on $W^{2,p}(0,1)\cap W^{1,p}_0(0,1)$. Furthermore $$\ker(\lambda_1 I - A_p) = \{r\sin(\pi (\cdot)) \ | \  r\in\R\} \quad\text{and}\quad f_\pm(x) = \pm \pi/2 \quad\text{for}\quad x\in (0,1).$$ Hence it is not difficult to verify that condition $(LL1)$ is holds for $\bar u\in X_0\setminus\{0\}$, where $X_0 = \ker(\lambda_1 I - A_p)$. On the other hand, writing $f(s) := -\arctan(s)$ for $s\in\R$, one can easily check that condition $(LL2)$ is also satisfied for $\bar u\in X_0\setminus\{0\}$, where $X_0 = \ker(\lambda_1 I - A_p)$. \hfill $\square$
\end{example}

The following theorem shows that the conditions $(G1)$ and $(G2)$ are also implicated by \emph{the strong resonance conditions}, studied for example in \cite{MR713209}, \cite{MR597281}.

\begin{theorem}\label{lem-est3}
Assume that there is a continuous function $f_\infty \colon \overline{\Omega} \to \mathbb{R}$ such that
\begin{equation}\label{asu}
f_\infty(x)  = \lim_{|s| \to +\infty} f(x,s)\cdot s \ffor x\in\Omega.
\end{equation}
\makebox[6mm][l]{(i)}\parbox[t]{117mm}{Assume that \\[-9pt]
\begin{equation*}\leqno{(SR1)}
\quad\left\{\begin{aligned}
\parbox[t][][t]{102mm}{\em there is a function $h\in L^1(\Omega)$ such that $$\hspace{-80pt}f(x,s)\cdot s \ge h(x) \ \  \text{ for } \ \  (x,s)\in \Omega\times\R  \ \text{ and }  \ \int_\Omega f_\infty(x)\d x > 0.$$}
\end{aligned} \right.
\end{equation*}
If the sets $B_1\subset X^\alpha_+\oplus X^\alpha_-$ and $B_2 \subset X_0$ are bounded in the norms $\|\cdot\|_\alpha$ and $\|\cdot\|_{L^2}$, respectively, then there are constants $R,\rho > 0$ such that
$$\<F(\z + \x), \x\>_{L^2} > - \<F(\z + \x), \y\>_{L^2} + \rho$$
for any $(\y, \z,\x)\in B_1\times B_2 \times X_0$, with $\|\x\|_{L^2}\ge R$.}\\[5pt]
\makebox[6mm][l]{(ii)}\parbox[t]{117mm}{Assume that \\[-9pt]
\begin{equation*}\leqno{(SR2)}
\quad\left\{\begin{aligned}
\parbox[t][][t]{102mm}{\em there is a function $h\in L^1(\Omega)$ such that $$\hspace{-80pt}f(x,s)\cdot s \le h(x) \ \  \text{ for } \ \  (x,s)\in \Omega\times\R  \ \text{ and }  \ \int_\Omega f_\infty(x)\d x < 0.$$}
\end{aligned} \right.
\end{equation*}
If the sets $B_1\subset X^\alpha_+\oplus X^\alpha_-$ and $B_2 \subset X_0$ are bounded in the norms $\|\cdot\|_\alpha$ and $\|\cdot\|_{L^2}$, respectively, then there are constants $R,\rho > 0$ such that
$$\<F(\z + \x), \x\>_{L^2} < - \<F(\z + \x), \y\>_{L^2} - \rho$$
for any $(\y, \z,\x)\in B_1\times B_2 \times X_0$, with $\|\x\|_{L^2}\ge R$.}\\[5pt]
\end{theorem}

\begin{remark}
Let us observe that under the hypotheses of Theorem \ref{lem-est3}, one has
\begin{equation*}
f_\pm (x) = \lim_{s \to \pm\infty} f(x,s) = 0 \text{ \ for \ } x\in\Omega,
\end{equation*}
which implies that the Landesman-Lazer conditions $(LL1)$ and $(LL2)$ used in Theorem \ref{lem-est2} are not valid in this case. \hfill $\square$
\end{remark}
\noindent\textbf{Proof of Theorem \ref{lem-est3}.} Similarly as before, with some technical modifications, the argument goes in the lines of the proof of \cite[Theorem 6.9]{Kok2}. We encourage the reader one more time to reconstruct details. \hfill $\square$ \\[-5pt]

\begin{example}
We describe an example of the situation where the strong resonance conditions $(SR1)$ and $(SR2)$ are satisfied. \\[5pt]
Let us note that, if $f:\R\to\R$ is a map given by $f(s) := s/(1 + s^2)$ for $s\in \R$ then $f(s)\cdot s \to 1$ as $|s|\to +\infty$. Hence, one can easily check that condition $(SR1)$ is satisfied with $f_\infty \equiv 1$. On the other hand, writing $f(s) := - s/(1 + s^2)$ for $s\in \R$, one has $f(s)\cdot s \to -1$ as $|s|\to +\infty$, and consequently condition $(SR2)$ is satisfied with $f_\infty \equiv - 1$. \hfill $\square$
\end{example}

\subsection{Criteria on existence of maximal compact invariant sets}

Let us observe that the equation \eqref{row-abs} can by written as
\begin{equation}\label{row-dr2}
\dot w(t) = -\A_p w(t) + \F (w(t)),  \qquad   t > 0.
\end{equation}
where $\A_p:\E\supset D(\A_p)\to\E$ is a linear operator on $\E:= X^\alpha\times X$ given by
\begin{equation*}
\begin{aligned}
D(\A_p) & :=\{(\x,\y)\in X^\alpha\times X \ | \ \x + c \y\in D(A_p)\} \\
\A_p(\x,\y) & :=(-\y, A_p(\x + c \y) - \lambda \x)
\end{aligned}
\end{equation*}
and $\F: \E\to\E$ is a map defined by $\F(\x,\y):= (0, F(\x))$ for $(\x,\y)\in \E$. Then Remark \ref{rem-pom2} $(b)$ and Proposition \ref{th-exist1} $(a)$ assert that the semiflow $\Pphi:[0,+\infty)\times\E\to\E$ associated with the equation \eqref{row-dr2} is well-defined, continuous and any bounded subset of $\E$
is admissible with respect to $\Pphi$. Let us first prove the following {\em criterion with Landesman-Lazer type conditions}.
\begin{theorem}\label{th-crit-ogrh}
Let $f_+,f_-\colon \Omega \to \mathbb{R}$ be continuous functions such that
\begin{equation*}
f_+(x) = \lim_{s \to +\infty} f(x,s) \quad\text{and}\quad f_-(x) = \lim_{s \to -\infty} f(x,s) \quad\text{for}\quad x\in\Omega.
\end{equation*}
Assume that $\lambda = \lambda_k$ is the $k$-th eigenvalue of the operator $A_p$ and let $K_\infty$ be the union of all bounded orbits of the semiflow $\Pphi$. Then, there is a closed isolating neighborhood $N\subset \E$, admissible with respect to $\Pphi$, such that $K_\infty=\inv(N,\Pphi)$ and the following assertions hold:
\begin{enumerate}
\item[(i)] if condition $(LL1)$ is satisfied, then $h(\Pphi, K_\infty) = \Sigma^{d_k}$, \\[-9pt]
\item[(ii)] if condition $(LL2)$ is satisfied, then $h(\Pphi, K_\infty) = \Sigma^{d_{k-1}}$.
\end{enumerate}
Here $d_l:=\sum_{i=1}^l \dim\ker(\lambda_i I - A_p)$ for $l\ge 1$ with the exceptional case $d_0:= 0$.
\end{theorem}
\noindent\textbf{Proof.} For the proof it is enough to combine Theorem \ref{lem-est2} and Theorem \ref{th-lan-laz}. \hfill $\square$ \\

Now we proceed to the following \emph{criterion with strong resonance conditions}.

\begin{theorem}\label{th-crit-ogrsrh}
Let $f_\infty \colon \o\Omega \to \mathbb{R}$ be a continuous function such that
\begin{gather*}
f_\infty(x)  = \lim_{|s| \to +\infty} f(x,s)\cdot s \ffor x\in\Omega.
\end{gather*}
Assume that $\lambda = \lambda_k$ is the $k$-th eigenvalue of the operator $A_p$ and let $K_\infty$ be the union of all bounded orbits of the semiflow $\Pphi$. Then, there is a closed isolating neighborhood $N\subset \E$, admissible with respect to $\Pphi$, such that $K_\infty=\inv(N,\Pphi)$ and the following assertions hold:
\begin{enumerate}
\item[(i)] if condition $(SR1)$ is satisfied, then $h(\Pphi, K_\infty) = \Sigma^{d_k}$, \\[-9pt]
\item[(ii)] if condition $(SR2)$ is satisfied, then $h(\Pphi, K_\infty) = \Sigma^{d_{k-1}}$.
\end{enumerate}
Here $d_l:=\sum_{i=1}^l \dim\ker(\lambda_i I - A_p)$ for $l\ge 1$ with the exceptional case $d_0:= 0$.
\end{theorem}
\noindent\textbf{Proof.} For the proof it is enough to combine Theorem \ref{lem-est3} and Theorem \ref{th-lan-laz}. \hfill $\square$

\begin{remark}
By Theorems \ref{th-crit-ogrh} and \ref{th-crit-ogrsrh} and existence property of Conley index, it follows that the set $K_\infty$ for the semiflow $\Pphi$ is compact and nonempty provided either Landesman-Lazer or strong resonance conditions are satisfied.
\end{remark}

\subsection{Criteria on existence of connecting orbits}

We continue our studies on the equation \eqref{A-eps-res-ah}. Here we assume that the map $f:\o\Omega\times \R\to\R$ is of class $C^1$ and, in addition do the conditions $(E1)$ and $(E2)$, we require that \\[4pt]
\noindent\makebox[22pt][l]{$(E3)$} \parbox[t][][t]{115mm}{$f(x,0) = 0$ for $x\in\Omega$ and there is $\nu\in\R$ such that $\nu = D_s f(x,0)$ for $x\in\Omega$.}\\[5pt]
\indent From assumption $(E3)$ it follows that $\Pphi(t,0) = 0$ for $t\ge 0$. The following criterion with Landesman-Lazer conditions determines when the maximal compact invariant set $K_\infty$ obtained in Theorem \ref{th-crit-ogrh} contains an orbit which is either homoclinic at zero or connects zero with another invariant set.
\begin{theorem}
Assume that there are continuous maps $f_+,f_-\colon \Omega \to \mathbb{R}$ such that
\begin{equation*}
f_+(x) = \lim_{s \to +\infty} f(x,s) \quad\text{and}\quad f_-(x) = \lim_{s \to -\infty} f(x,s) \quad\text{for}\quad x\in\Omega.
\end{equation*}
If $\lambda=\lambda_k$ is the $k$-th eigenvalue of the operator $A_p$, then there is a non-zero solution $w:\R \to \E$ of the semiflow $\Pphi$ such that $$w(\R)\subset K_\infty \quad\text{and either}\quad \lim_{t\to -\infty} w(t) = 0 \quad\text{or}\quad \lim_{t\to +\infty} w(t) = 0,$$ provided one of the following statements is satisfied: \\[-11pt]
\begin{enumerate}
\item[(i)] condition $(LL1)$ holds and $\lambda_l < \lambda + \nu < \lambda_{l+1}$ where $\lambda_l \neq \lambda$,\\[-9pt]
\item[(ii)] condition $(LL1)$ holds and $\lambda + \nu < \lambda_1$, \\[-9pt]
\item[(iii)] condition $(LL2)$ holds, $\lambda_{l-1} < \lambda + \nu < \lambda_l$ and $\lambda \neq \lambda_l$, where $l\ge 2$, \\[-9pt]
\item[(iv)] condition $(LL2)$ holds, $\lambda + \nu < \lambda_1$ and $\lambda \neq \lambda_1$.
\end{enumerate}
\end{theorem}

\noindent\textbf{Proof.} From assumption $(E3)$ one can easily prove that $\F$ is differentiable at $0$ and its derivative $D\F(0)\in L(\E, \E)$ has the following form $$D\F(0)[\x,\y] = \nu (0,\x) \ffor (\x,\y)\in \E.$$
Therefore, if $\lambda + \nu \notin\sigma(A_p)$, then \cite[Theorem 3.5]{MR910097} implies that $\{0\}$ is isolated invariant set and $h(\Pphi,\{0\})= \Sigma^{b_l}$, where $b_l := 0$ if $\lambda + \nu < \lambda_1$ and $$b_l := \sum_{i=1}^l \dim\ker(\lambda_i I - A) \ \ \text{ if } \ \ \lambda_l < \lambda + \nu < \lambda_{l+1}.$$
Combining this with Theorems \ref{th-lan-laz}, \ref{lem-est2} and \cite[Theorem 11.5]{MR910097} gives desired assertion. \hfill $\square$ \\

The following criterion is similar to the previous one, with the difference that strong resonance conditions are used to determine the existence of connecting orbit contained in the set $K_\infty$.

\begin{theorem}
Assume that there is a continuous function $f_\infty \colon \o\Omega \to \mathbb{R}$ such that
\begin{gather*}
f_\infty(x)  = \lim_{|s| \to +\infty} f(x,s)\cdot s \quad\text{for} \quad x\in\Omega.
\end{gather*}
If $\lambda=\lambda_k$ is the $k$-th eigenvalue of the operator $A_p$, then there is a non-zero solution $w:\R \to \E$ of the semiflow $\Pphi$ such that $$w(\R)\subset K_\infty \quad\text{and either}\quad \lim_{t\to -\infty} w(t) = 0 \quad\text{or}\quad \lim_{t\to +\infty} w(t) = 0,$$ provided one of the following statements is satisfied: \\[-11pt]
\begin{enumerate}
\item[(i)] condition $(SR1)$ holds and $\lambda_l < \lambda + \nu < \lambda_{l+1}$ where $\lambda_l \neq \lambda$; \\[-9pt]
\item[(ii)] condition $(SR1)$ holds and $\lambda + \nu < \lambda_1$; \\[-9pt]
\item[(iii)] condition $(SR2)$ holds and $\lambda_{l-1} < \lambda + \nu < \lambda_l$ where $\lambda \neq \lambda_l$, $l\ge 2$; \\[-9pt]
\item[(iv)] condition $(SR2)$ holds and $\lambda + \nu < \lambda_1$. \\[-9pt]
\end{enumerate}
\end{theorem}

\def\cprime{$'$} \def\polhk#1{\setbox0=\hbox{#1}{\ooalign{\hidewidth
  \lower1.5ex\hbox{`}\hidewidth\crcr\unhbox0}}} \def\cprime{$'$}
  \def\cprime{$'$} \def\cprime{$'$}
\providecommand{\bysame}{\leavevmode\hbox to3em{\hrulefill}\thinspace}
\providecommand{\MR}{\relax\ifhmode\unskip\space\fi MR }
\providecommand{\MRhref}[2]{%
  \href{http://www.ams.org/mathscinet-getitem?mr=#1}{#2}
}
\providecommand{\href}[2]{#2}

\parindent = 0 pt

\end{document}